\documentclass[10pt]{article}

\usepackage{amsmath}
\usepackage{amssymb}
\usepackage{indentfirst}
\usepackage{graphics} 
\usepackage{color}

\makeatletter

\@addtoreset{equation}{section}
\makeatother
\def\<{\langle}
\def\>{\rangle}

\newtheorem{lem}{Lemma}[section]
\newtheorem{theo}{Theorem}[section]
\newtheorem{rem}{Remark}[section]
\newtheorem{pro}{Proposition}[section]

\makeatletter
   
   \@addtoreset{equation}{section}
\makeatother

\begin{document}
\title{\bf  A dissipative logarithmic type evolution equation:\\ asymptotic profile and optimal estimates}

\author{ Alessandra Piske\thanks {alessandrapiske@gmail.com} \; and \;  Ruy Coimbra Char\~ao\thanks{Corresponding author: ruy.charao@ufsc.br}   \\{\small Department of Mathematics} \\
\small Graduate Program on Pure and Applied Mathematics\\
{\small Federal University of Santa Catarina} \\ {\small 88040-270, Florianopolis, Brazil,} 
\\
and\\Ryo Ikehata\thanks{ikehatar@hiroshima-u.ac.jp} \\ {\small Department of Mathematics}\\ {\small Division of Educational Sciences}\\ {\small Graduate School of Humanities and Social Sciences} \\ {\small Hiroshima University} \\ {\small Higashi-Hiroshima 739-8524, Japan}}
\date{}
\maketitle
\begin{abstract}
We introduce a new model of the  logarithmic type of  wave like equation with a nonlocal logarithmic damping mechanism, which is rather weakly effective as compared with frequently studied fractional damping cases.  We consider the Cauchy problem for this new model in ${\bf R}^{n}$, and study the asymptotic profile and optimal decay and/or blowup rates of solutions as $t \to \infty$ in $L^{2}$-sense. The operator $L$ considered in this paper was used to dissipate the solutions of the wave equation  in the paper  studied by Char\~ao-Ikehata \cite{Log-damping} and in the low frequency parameters the principal part of the equation and the damping term is rather weakly effective than that of well-studied power type one such as $(-\Delta)^{\theta}u_{t}$ with $\theta \in (0,1]$.
\end{abstract}
\section{Introduction}
\footnote[0]{Keywords and Phrases: Wave-like equation; Logarithmic damping; $L^{2}$-decay; asymptotic profile, optimal estimates.}
\footnote[0]{2010 Mathematics Subject Classification. Primary 35L05; Secondary 35B40, 35C20, 35S05.}
We present and consider a new model of evolution equation with a logarithmic damping term:
\begin{align}
& u_{tt} + Lu + Lu_{t} = 0,\ \ \ (t,x)\in (0,\infty)\times {\bf R}^{n},\label{eqn}\\
& u(0,x)= u_{0}(x),\ \ u_{t}(0,x)= u_{1}(x),\ \ \ x\in{\bf R}^{n} ,\label{initial}
\end{align}
where the linear operator  
\[L: D(L) \subset L^{2}({\bf R}^{n}) \to L^{2}({\bf R}^{n})\]
is defined as follows: 
\[D(L) := \left\{f \in L^{2}({\bf R}^{n}) \,\bigm|\,\int_{{\bf R}^{n}}(\log(1+\vert\xi\vert^{2}))^{2}\vert\hat{f}(\xi)\vert^{2}d\xi < +\infty\right\},\]
\[(Lf) (x) := {\cal F}_{\xi\to x}^{-1}\left(\log (1+\vert\xi\vert^{2})\hat{f}(\xi)\right)(x), \quad \; f \in D(L). \]
\noindent
Here, one has just denoted the Fourier transform ${\cal F}_{x\to\xi}(f)(\xi)$ of $f(x)$ by 
\[{\cal F}_{x\to\xi}(f)(\xi) = \hat{f}(\xi) := \displaystyle{\int_{{\bf R}^{n}}}e^{-ix\cdot\xi}f(x)dx, \quad \xi \in {\bf R}^n,\]
as usual with $i := \sqrt{-1}$, and ${\cal F}_{\xi\to x}^{-1}$ expresses its inverse Fourier transform. 
\noindent
Since the operator $L$ is non-negative and self-adjoint in $L^{2}({\bf R}^{n})$ (see \cite{Log-damping}), the square root 
$$L^{1/2}: D(L^{1/2}) \subset L^{2}({\bf R}^{n}) \to L^{2}({\bf R}^{n})$$
can be defined, and is also nonnegative and self-adjoint with its domain
\[D(L^{1/2}) = \left\{f \in L^{2}({\bf R}^{n}) \,\bigm|\,\int_{{\bf R}^{n}}\log(1+\vert\xi\vert^{2})\vert\hat{f}(\xi)\vert^{2}d\xi < +\infty\right\}.\] 
Note that $D(L^{1/2})$ becomes Hilbert space with its graph norm
$$\Vert v\Vert_{D(L^{1/2})} := \left(\Vert v\Vert^{2}_{L^{2}} + \Vert L^{1/2}v\Vert^{2}_{L^{2}}\right)^{1/2}.$$
It is easy to check that 
$$H^{s}({\bf R}^{n}) \hookrightarrow D(L^{1/2}) \hookrightarrow  L^{2}({\bf R}^{n})$$
for $s > 0$.\\
\noindent
Symbolically writing, one can see
\[L = \log(I-\Delta),\]
where $\Delta$ is the usual Laplace operator defined on $H^2({\bf R}^n)$.

Now, for the time being, we choose the initial data $(u_{0},u_{1})$ as follows:
\[u_{0} \in D(L^{1/2}),\quad u_{1} \in L^{2}({\bf R}^{n}).\]
The existence of the unique solution to problem \eqref{eqn}-\eqref{initial} can be discussed by employing a similar argument to \cite[Proposition 2.1]{ITY} based on Lumer-Phillips Theorem, and one can find that the problem \eqref{eqn}-\eqref{initial} has a unique mild solution
\[u \in C([0,\infty);D(L^{1/2})) \cap C^{1}([0,\infty);L^{2}({\bf R}^{n}))\]
satisfying the energy inequality
\begin{equation}\label{energy}
E_{u}(t) \leq E_{u}(0),
\end{equation} 
where
\[
E_{u}(t) := \frac{1}{2}\left(\Vert u_{t}(t,\cdot)\Vert_{L^{2}}^{2} + \Vert \log^{1/2}(I-\Delta) u(t,\cdot)\Vert_{L^{2}}^{2}\right).
\]
The inequality \eqref{energy} implies that the the total energy is a non increasing function in time because of the existence of some kind of dissipative term $Lu_{t}$.  

A main topic of this paper is to find an asymptotic profile of solutions in the $L^{2}$ framework to problem  \eqref{eqn}-\eqref{initial}, and to apply it to investigate the optimal rate of decay of solutions in terms of the $L^{2}$-norm. We study the equation \eqref{eqn} only from the pure mathematical point of view.\\  

A motivation of this research has its origin in the study of the strongly damped wave equation:
\begin{equation}\label{strong}
u_{tt}-\Delta u -\Delta u_{t} = 0.
\end{equation}
An analysis of the dissipative mechanism of \eqref{strong} goes back to the two pioneering works of G. Ponce \cite{Po} and Y. Shibata \cite{S}, where they studied various $L^{p}$-$L^{q}$ estimates of the solution to the Cauchy problem of \eqref{strong}. After them, an asymptotic profile and the optimal estimates of the solution can be introduced in the papers \cite{ITY}, \cite{I-14} and \cite{IO}. They investigated a singularity near $0$-frequency region of the solution to \eqref{strong} in terms of $L^{2}$-norm of solutions. In this connection, in \cite{B, BV-1, BV-2} and \cite{Mi} a higher order asymptotic expansion of the solution as $t \to \infty$ to the equation \eqref{strong} is precisely investigated.

On the other hand, the so-called critical exponent problem for semi-linear equations of \eqref{strong} is first developed by D'Abbicco-Reissig \cite{DR}, and this paper has been the beginning of a series of related papers studying structurally damped wave models with nonlinearity. Unfortunately, at present nobody knows the precise value of the critical exponent $p^{*}$ of the equation \eqref{strong} with power type nonlinearity $\vert u\vert^{p}$. A study in \cite{DR} is based on the $L^{p}-L^{q}$-estimates derived in \cite{S}.

Recently, the equation \eqref{strong} is generalized to the linear and semi-linear models, respectively:
\begin{equation}\label{strong-2}
u_{tt} + (-\Delta)^{\sigma}u + (-\Delta)^{\theta}u_{t} = 0,
\end{equation}
\begin{equation}\label{strong-3}
u_{tt} + (-\Delta)^{\sigma}u + (-\Delta)^{\theta}u_{t} = f(u,u_{t}).
\end{equation}
A study on asymptotic profile and $L^{p}$-$L^{q}$ estimates to the equation \eqref{strong-2} has been done in the papers \cite{CLI}, \cite{DE}, \cite{DEP}, \cite{DGL}, \cite{LIC}, \cite{NR}, and \cite{IT}, and the corresponding critical exponent problems (mainly) to the equation \eqref{strong-3} are treated in the papers \cite{AR, AR-2}, \cite{DE-2}, \cite{FC}, \cite{K}, and \cite{PMR}. 

In \cite{II} and \cite{FIM}, the so-called regularity-loss structure of the solution in the high frequency zone can be studied to the equation \eqref{strong-2} with $\sigma = 1$ and $\theta > 1$, and these researches are strongly inspired from the abstract theory due to \cite{GGH}. Such a regularity-loss structure has been first discovered by S. Kawashima through the analysis for dissipative Timoshenko system. A more general model than \eqref{strong-2} is studied by \cite{C-E-I}. The aim of that work in \cite{C-E-I} is to obtain asymptotic profile and optimal decay rates in case of a super damping (i.e., $\sigma < \theta$). We have much more interesting results about more generalized evolution equations such as memory type of damping, double one, rotational inertia term case, and etc... than \eqref{strong-3}, however, we do not mention them not to spread in vain our topics.

Quite recently Char\~ao-Ikehata  \cite{Log-damping} introduced a new type of damping term of logarithm type to the wave equation, and it is expressed in the Fourier space as follows:
\begin{equation}\label{logdamping}
\hat{u}_{tt} + \vert\xi\vert^{2}\hat{u} + \log(1+\vert\xi\vert^{2})\hat{u}_{t} = 0.
\end{equation}
Symbolically writing, one sees
\[u_{tt}-\Delta u + \log(I-\Delta)u_{t} = 0.\]
In \cite{CDI}, \eqref{logdamping} is more generalized to the equation such that
\[u_{tt}-\Delta u + \log(I + (-\Delta)^{\theta})u_{t} = 0\]
for $\theta > 1/2$. As is easily seen that the characteristic roots $\lambda_{\pm}$ for the characteristic polynomial of \eqref{logdamping} such that 
\[\lambda^{2} + \log(1+\vert\xi\vert^{2})\lambda + \vert\xi\vert^{2} = 0\]
are complex-valued for all $\xi \in {\bf R}^{n}$, although the contribution on the decay structure from the high frequency parameters is very small.  In this connection, in \cite{WS} they study another model with double dispersion for which oscillations appear at both low and high frequencies. 

On reconsidering our problem \eqref{eqn}-\eqref{initial} in the Fourier space, our equation becomes  
\begin{equation}\label{log-logdamping}
\hat{u}_{tt} + \log(1+\vert\xi\vert^{2})\hat{u} + \log(1+\vert\xi\vert^{2})\hat{u}_{t} = 0.
\end{equation} 
Characteristics roots of \eqref{log-logdamping} are complex-valued only for small $\xi \in {\bf R}^{n}$, and in the large frequency zone, the roots are real-valued, and this is similar to the strong damping case \eqref{strong}.

To get started, we first investigate the decay rate of the total energy $E_{u}(t)$ and $L^{2}$-norm of the solution itself under the $L^{1}$-framework on the initial data.
\begin{pro}\label{L-energy}
Let $u(t,x)$ be the solution to problem \eqref{eqn}-\eqref{initial} with initial data
$$(u_0, u_1) \in \left(D(L^{1/2})\cap L^1({\bf R}^n)\right)\times \left(L^2({\bf R}^n)\cap L^1({\bf R}^n)\right).$$
Then, the total energy of this system satisfies for $t \gg 0$
$$ \left \| u_t (t, \cdot) \right \|_{L^2}^{2} + \left \| L^{1/2}u (t, \cdot) \right \|_{L^2}^{2} 
 \leq C_n \left(\left \| u_1 \right \|_{L^{1}}^{2}t^{- \frac{n}{2}} + \left \| u_0 \right \|_{L^{1}}^{2}  t^{- \frac{n+2}{2}} \right) +
2^{-\frac{t}{4}} \left( \left \| u_1 \right \|_{L^{2}}^{2} +  \left \|u_0 \right \|_{L^{2}}^{2}\right) +2e^{-\frac{t}{4}}E_u(0).
$$
\end{pro}
\begin{rem}
{\rm The above proposition says that the total energy of the system decays as $t^{-n/2}$, that is
 $$ E_u(t) \leq 
 C_{1,n}  \left( E_u(0) + \left \| u_0 \right\|_{L^{2}}^{2} +  \left \| u_0 \right \|_{L^{1}}^{2} +  \left \| u_1 \right \|_{L^{1}}^{2} \right )t^{- \frac{n}{2}} ,  \quad t \gg 1,$$
 with a constant $C_{1,n}>0$ depending only on $n$.}
 \end{rem}
 \begin{pro}\label{L2-u} Let $n >2$ and $u(t,x)$ be the solution to problem \eqref{eqn}-\eqref{initial} with initial data 
$$  u_0,u_1 \in L^2({\bf R}^n)\cap L^1({\bf R}^n).$$
Then
$$   \left \| u (t, \cdot) \right \|_{L^{2}} \leq C_n  \left (    \left \| u_0 \right \|_{L^{2}} +  \left \| u_1 \right \|_{L^{2}}+  \left \| u_0 \right \|_{L^{1}} +  \left \| u_1 \right \|_{L^{1}} \right )t^{- \frac{n-2}{4}}, \quad  t \gg 1, $$
with a constant $C_n>0$ depending only on $n$.
\end{pro} 
\begin{rem}\,{\rm The decay rate of the quantity $\Vert u(t,\cdot)\Vert_{L^{2}}$ can be derived only for the spatial dimension $n > 2$ under the $L^{1}$-regularity on the initial data. $n = 1,2$ cases have a strong singularity near $0$-frequency region. This singularity can be observed in the following main results below.}
\end{rem} 
 
In order to investigate the optimality of decay rates of the quantity $\Vert u(t,\cdot)\Vert_{L^{2}}$ just obtained in Proposition \ref{L2-u} we do study the asymptotic profile of the solution $u(t,x)$ as $t \to \infty$ in $L^{2}$-sense. Our new result reads as follows. At this stage, it suffices to assume $u_{0} = 0$ without loss of generality. 

\begin{theo}\label{main-theo}
\, Let $n \geq 1$, and let $u_{0} = 0$, and $u_{1} \in \left(L^{2}({\bf R}^{n})\cap L^{1,1}({\bf R}^{n})\right)$. Then, the unique solution $u(t,x)$ to problem  \eqref{eqn}-\eqref{initial} satisfies
\[\left\Vert u(t,\cdot) - \left(\int_{{\bf R}^{n}}u_{1}(x)dx\right){\cal F}_{\xi\to x}^{-1}\left((1+\vert\xi\vert^{2})^{-\frac{t}{2}}\frac{\sin(\sqrt{\log(1+|\xi|^2)}\;t))}{\sqrt{\log(1+\vert\xi\vert^2)}}\right)\right\Vert_{L^{2}} \leq I_{0}t^{-\frac{n}{4}}, \quad (t \gg 1),\]
where
\[I_{0} := \Vert u_{1}\Vert_{L^{2}} + \Vert (1+\vert x\vert)u_{1}\Vert_{L^{1}}.\]
\end{theo}
\begin{rem}\,{\rm If we apply the general theory developed in \cite{ITY} to the abstract evolution equation 
\begin{equation}\label{abstract}
u_{tt} + Au + Au_{t} = 0,
\end{equation}
where $A$ is a nonnegative self-adjoint operator in a real Hilbert space, at least one can observe the asymptotic profile of the solution to problem \eqref{eqn}-\eqref{initial} is
\[e^{-tL/2}\frac{\sin(L^{1/2}t)}{L^{1/2}}u_{1}.\]
By restricting the initial data further to the class $L^{1,1}({\bf R}^{n})$, one can obtain the statement of Theorem 1.1.}
\end{rem}

By using Theorem \ref{main-theo} the optimal estimates in $t$ of the solution $u(t,x)$ can be derived in terms of $L^{2}$-norm. As a gain of Theorem \ref{main-theo} one can treat $n \geq 3$ and $n = 1,2$ as well. Note that $n = 1,2$ are missing dimensions in Proposition \ref{L2-u}. We set
\[P_{1} := \int_{{\bf R}^{n}}u_{1}(x)dx.\]
The most newest and important task of this paper is to get the asymptotic behavior in $t$ of the following improper integrals such that 
\begin{equation}\label{improper}
{\cal I}_{n}(t) := \int _{{\bf R}^n} \frac{(1+|\xi |^2)^{-t} \sin ^2 (t\sqrt{ \log (1+|\xi |^2)})}{\log (1+|\xi |^2)} d \xi. 
\end{equation}
This type of improper integrals have its origin in Theorem \ref{main-theo} and the problem to capture precise decay and/or blowup orders seem to be purely new from a mathematical point of view, and it seems extremely interesting to deal with the low dimensional case $n = 1,2$ because the integrals \eqref{improper} include a log-type singularity near $\xi = 0$. Our important contribution is as follows. In the $1$ and $2$ dimensional cases one has strong singularities near $\xi = 0$ even in the log-damping case. 
\begin{theo}\label{main-theo2}
\, Let $n \geq 1$, and let $u_{0} = 0$, and $u_{1} \in \left(L^{2}({\bf R}^{n})\cap L^{1,1}({\bf R}^{n})\right)$. Then, the unique solution $u(t,x)$ to problem {\rm (1.1)}-{\rm (1.2)} satisfies the following properties:\\
\vspace{0.1cm}
{\rm (i)}\,if $n = 1$, then $C_{1}\vert P_{1}\vert \sqrt{t} \leq \Vert u(t,\cdot)\Vert_{L^{2}} \leq C_{1}^{-1}I_{0}\sqrt{t}$ {\rm (}$t \gg 1${\rm )},\\
\vspace{0.1cm}
{\rm (ii)}\,\,if $n = 2$, then $C_{2}\vert P_{1}\vert \sqrt{\log t} \leq \Vert u(t,\cdot)\Vert_{L^{2}} \leq C_{2}^{-1}I_{0}\sqrt{\log t}$ {\rm (}$t \gg 1${\rm )},\\
\vspace{0.1cm}
{\rm (iii)}\,\,\,if $n \geq 3$, then $C_{n}\vert P_{1}\vert t^{-\frac{n-2}{4}} \leq \Vert u(t,\cdot)\Vert_{L^{2}} \leq C_{n}^{-1}I_{0}t^{-\frac{n-2}{4}}$ {\rm (}$t \gg 1${\rm )}.\\
Here $I_{0}$ is a constant defined in Theorem {\rm \ref{main-theo}}, and $C_{n}$ {\rm (}$n \in {\bf N}${\rm )} are constants independent from any $t$ and initial data.  
\end{theo}
\begin{rem}\,{\rm As a result, all estimates derived in Theorem \ref{main-theo2} are overlapped already known results in \cite{IO} and/or \cite{Log-damping}, and this is quite natural because $\log (1+\vert\xi\vert^{2}) \sim \vert\xi\vert^{2}$ for small $\xi \in {\bf R}^{n}$, and the main contribution to the estimates above comes from the low frequency region in $\xi \in {\bf R}^{n}$. However, by replacing the operator $A = -\Delta$ to $L = \log(I-\Delta)$ in the equation \eqref{abstract}, we encounter a big obstacle when one gets such estimates stated in Theorem \ref{main-theo2}, and this difficulty comes from the way that how we treat the improper integral \eqref{improper}. A big technical difficulties occur.  }
\end{rem}

\par
\vspace{0.1cm}

This paper is organized as follows. In section 2 we prepare several important lemmas, which will be used later, and in particular, these lemmas are closely related with hypergeometric functions (see \cite{Log-damping}). Propositions \ref{L-energy} and \ref{L2-u} can be proved in Section 3 based on the energy method due to \cite{UKS}. In Section 4, we derive the leading term (as $t \to \infty$) of the solution to problem \eqref{eqn}-\eqref{initial}. Section 5 is devoted to the derivation of the optimal decay rate of the $L^{2}$-norm of the solution in case of $n \geq 3$, and $1$ and $2$ dimensional cases for the optimality of the $L^{2}$-estimates of the solution will be investigated in Sections 6.\\

{\bf Notation.} {\small Throughout this paper, $\| \cdot\|_q$ stands for the usual $L^q({\bf R}^{n})$-norm. For simplicity of notation, in particular, we use $\| \cdot\|$ instead of $\| \cdot\|_2$. Furthermore, we denote $\Vert\cdot\Vert_{H^{l}}$ as the usual $H^{l}$-norm. Furthermore, we define a relation $f(t) \sim g(t)$ as $t \to \infty$ by: there exist constant $C_{j} > 0$ ($j = 1,2$) such that
\[C_{1}g(t) \leq f(t) \leq C_{2}g(t)\quad (t \gg 1).\] 

We also introduce the following weighted functional spaces.
\[L^{1,\gamma}({\bf R}^{n}) := \left\{f \in L^{1}({\bf R}^{n}) \; \bigm| \; \Vert f\Vert_{1,\gamma} := \int_{{\bf R}^{n}}(1+\vert x\vert^{\gamma})\vert f(x)\vert dx < +\infty\right\}.\]
Finally, we denote the surface area of the $n$-dimensional unit ball by $\omega_{n} := \displaystyle{\int_{\vert\omega\vert = 1}}d\omega$. 

}


\section{General basic results}
In this section we shall prepare important lemmas to derive precise estimates of the various quantities related to the solution to problem \eqref{eqn}-\eqref{initial}. These are already studied and developed in our previous works (see \cite{Log-damping, CDI}).

The following main estimate for the function  
\begin{equation}\label{2.1}
I_p(t)= \int_0^{1}(1+r^{2})^{-t}r^p dr
\end{equation} 
is a direct consequence of the cases $p \geq 0$ in Char\~ao-Ikehata \cite{Log-damping} and $-1<p<0$ in Char\~ao-D'Abbicco-Ikehata \cite{CDI}.

\begin{lem}\label{general-p}
{\rm Let $p > -1$ be a real number. Then
$$I_p(t) \sim t^{-\frac{p+1}{2}}, \quad t \gg 1.$$}
\end{lem}


In order to deal with the high frequency part of estimates, one relies on the function again 
\begin{equation}\label{2.2}
J_p(t)=\int_1^{\infty}(1+r^2)^{-t}r^p dr
\end{equation}
for $p \in {\bf R}$.

Then the next lemma is important to get estimates on the zone of high frequency to problem \eqref{eqn}--\eqref{initial}. The proof  appears in Char\~ao-Ikehata \cite{Log-damping}.
\begin{lem}\label{infit}
\,Let $p \in {\bf R}$. Then it holds that 
$$J_p(t) \sim \dfrac{2^{-t}}{t-1}, \quad t \gg 1.$$
\end{lem}
\vspace{0.2cm}
For later use we prepare the following simple lemma, which implies the exponential decay estimates of the middle frequency part.
\begin{lem}\label{intermid}\,Let $p \in {\bf R}$, and $\eta \in (0,1]$. Then there is a constant $C > 0$ such that 
$$\int_{\eta}^{1}(1+r^{2})^{-t}r^{p}dr \leq C(1+\eta^{2})^{-t}, \quad t \geq 0.$$
\end{lem}

\begin{rem}
{\rm We note that the proof of Lemma \ref{general-p} are proved using simple differential calculus and the theory from hypergeometric functions (see Watson \cite{W}).}

\end{rem}
\section{Asymptotic behavior via multiplier method}

In this section, we shall obtain optimal estimates of the total energy of the following Fourier transformed equation together with initial data of the original system \eqref{eqn}-\eqref{initial}. To do so we employ the so-called energy method in the Fourier space developed in \cite{UKS}. It seems to be a new development for this type of equation with logarithmic operators.
\begin{align}
&\widehat{u}_{tt} + \log (1+ | \xi |^2) \widehat{u} + \log (1+ | \xi |^2) \widehat{u}_t =0,  \quad  (t,\xi) \in (0, \infty) \times {\bf R}^n, \label{eqn2} \\ 
& \widehat{u}(0,\xi) =  \widehat{u}_0(\xi), \quad 
 \widehat{u}_t(0,\xi) =  \widehat{u}_1(\xi), \quad   \xi \in {\bf R}^n.  \label{initial2}
\end{align}

Multiplying the equation  \eqref{eqn2} by $ \bar{\widehat{u}_t} $ one can get the following point wise energy identity 
\begin{equation}\label{energy-id}
\frac{\mathrm{d} E_0(t,\xi) }{\mathrm{d} t} +  \log (1+ | \xi |^2) |\widehat{u}_t(t,\xi)|^{2} = 0, 
\end{equation}
where 
$$ E_0(t, \xi)= \frac{|\widehat{u}_t(t, \xi)|^2}{2} + \log (1+ | \xi |^2)\frac{|\widehat{u}(t, \xi)|^2}{2}, $$  
for $t >0$ and $\xi \in {\bf R}^n$, is the total  density of energy of the system \eqref{eqn2}-\eqref{initial2}. Note from \eqref{energy-id} that $ E_0(t,\xi)  $ is a decreasing function of $t$ for each $\xi$. 

Now we define the following function of $\xi$. The  way to choose the best $\rho(\xi)$-function is showed in the work by Luz-Ikehata-Char\~ao \cite{LIC}: 
\begin{equation}\label{ro}
 \rho (\xi) = \begin{cases}
\displaystyle{\frac{1}{2}} \log (1+|\xi |^2) & \text{ if }  | \xi | \leq \sqrt{e -1},  \\ 
\displaystyle{\frac{1}{2}}  & \text{ if } | \xi |  \geq \sqrt{e -1}.
\end{cases} 
\end{equation}

By multiplying the equation \eqref{eqn2} by $ \rho( \xi ) \bar{\widehat{u}} $ we obtain the identity
$$ \rho (\xi ) \frac{\mathrm{d}  }{\mathrm{d} t} \left ( \widehat{u}_t \bar{ \widehat{u}} \right ) - \rho (\xi ) |\widehat{u}_t|^2  + \log (1+ |\xi |^2 ) \rho(\xi) | \widehat{u}|^2 + \log (1+ |\xi |^2 ) \rho(\xi)  \frac{\mathrm{d}  }{\mathrm{d} t} \frac{|\widehat{u}|^2}{2} = 0,$$
for all $t>0$ and $\xi \in {\bf R}^n$. Taking the real part on the last identity we arrive at  
\begin{equation} \label{eq4}
\frac{\mathrm{d}  }{\mathrm{d} t} \left [ \rho (\xi ) \text{Re} \left ( \widehat{u}_t \bar{ \widehat{u}} \right ) + \rho(\xi) \log (1+ |\xi |^2 )   \frac{|\widehat{u}|^2}{2}   \right ] + \rho(\xi) \log (1+ |\xi |^2 )  | \widehat{u}|^2  =  \rho (\xi ) |\widehat{u}_t|^2 ,
\end{equation}
which holds for  $t>0$ and $\xi \in {\bf R}^n$.

To proceed further we need to define the following functions on $(0,\infty)\times {\bf R}^n$:
 \begin{align}\label{EFR}
& E(t, \xi) = E_0(t, \xi ) + \rho (\xi ) \text{Re} \left ( \widehat{u}_t(t, \xi ) \bar{ \widehat{u}}(t, \xi ) \right ) + \frac{\rho(\xi)}{2} \log (1+ |\xi |^2 )    |\widehat{u}(t, \xi )|^2, \nonumber
\\
& F(t, \xi ) =  \log (1+ | \xi |^2) |\widehat{u}_t(t, \xi )|^2 +  \rho(\xi) \log (1+ |\xi |^2 )  | \widehat{u}(t, \xi )|^2, 
\\
& R(t, \xi ) = \rho (\xi ) |\widehat{u}_t(t, \xi )|^2. \nonumber
\end{align}

Then, adding \eqref{energy-id} and \eqref{eq4}, we get the following  identity 
\begin{equation}\label{eq5}
 \frac{\mathrm{d} }{\mathrm{d} t} E(t, \xi) + F(t,\xi) = R(t, \xi), 
\end{equation}
which also holds for $t>0$ and $\xi \in {\bf R}^n$.
Before continuing our argument, we need the next remark.
\begin{lem}\label{rem-ro}
The function $\rho(\xi)$ defined in \eqref{ro} satisfies the estimates 
$$ \rho(\xi) \leq \frac{1}{2}  $$
for $|\xi| \leq \sqrt{e-1}$.
Moreover,  $$\rho^2(\xi) \leq \frac{1}{4} \log (1+|\xi |^2) $$
for all $\xi \in {\bf R}^n$.
\end{lem}
{\it Proof.}\,\,Indeed, for  $|\xi | \leq \sqrt{e -1} $ we have 
$\log (1+|\xi |^2)| \leq 1$ which implies 
$$\log (1+|\xi |^2) \leq \log^{\frac{1}{2}} (1+|\xi |^2).$$
Thus, $\rho ^2 ( \xi) \leq \frac{1}{4} \log (1+|\xi |^2)$
and $\rho (\xi ) \leq \frac{1}{2}$ according to the definition of $\rho(\xi)$ in \eqref{ro}. 

For   $|\xi | \geq \sqrt{e -1}$ one has  
$\log(1+|\xi |^2)  \geq 1$. Thus, $\rho ^2 (\xi) = \frac{1}{4}  \leq \frac{1}{4} \log (1+|\xi |^2)$.
\hfill
$\Box$


\begin{lem}\label{E-Eo}
$$ \frac{1}{2} E_0(t, \xi) \leq E(t, \xi) \leq 3 E_0(t, \xi), \quad t>0, \; \xi \in {\bf R}^n.$$
\end{lem}

{\it Proof.} Using the inequality   
$ \rho (\xi ) \text{Re} \left ( \widehat{u}_t \bar{ \widehat{u}} \right )  \geq - \frac{ |\widehat{u}_t |^2}{4} - \rho ^2(\xi) |\widehat{u}|^2 $  and Lemma \ref{rem-ro},
one has 
\begin{eqnarray*}
E(t, \xi) &=& E_0(t, \xi ) + \rho (\xi ) \text{Re} \left ( \widehat{u}_t \bar{ \widehat{u}} \right ) + \frac{\rho(\xi)}{2} \log (1+ |\xi |^2 )    |\widehat{u}|^2 \\
&\geq & E_0(t, \xi ) - \frac{ |\widehat{u}_t |^2}{4} - \rho^2(\xi) |\widehat{u}|^2 \\
&=& \frac{|\widehat{u}_t|^2}{2} + \log (1+ | \xi |^2)\frac{|\widehat{u}|^2}{2} - \frac{ |\widehat{u}_t |^2}{4} - \rho^2(\xi) |\widehat{u}|^2 \\
&=& \frac{1}{4} |\widehat{u}_t |^2 + \left (\frac{\log (1+ | \xi |^2)}{2}  - \rho^2(\xi)  \right )  |\widehat{u}|^2  \\
& \geq & \frac{1}{4} |\widehat{u}_t |^2 + \frac{\log (1+ | \xi |^2)}{4}  |\widehat{u}|^2  
=\frac{1}{2} E_0(t, \xi),
\end{eqnarray*}
which  holds for $t>0$ and $\xi \in {\bf R}^n$.

On the other hand, using Lemma \ref{rem-ro} one has the estimates 
\begin{eqnarray*}
E(t, \xi) &=& E_0(t, \xi ) + \rho (\xi ) \text{Re} \left ( \widehat{u}_t \bar{ \widehat{u}} \right ) + \frac{\rho(\xi)}{2} \log (1+ |\xi |^2 )    |\widehat{u}|^2 \\
& \leq & E_0(t, \xi ) + \frac{| \widehat{u}_t|^2}{2} + \frac{\rho^2(\xi)}{2} | \widehat{u}|^2 + \frac{\rho(\xi)}{2} \log (1+ |\xi |^2 )    |\widehat{u}|^2 \\
& \leq & E_0(t, \xi ) + \frac{| \widehat{u}_t|^2}{2} + \frac{ \log (1+|\xi |^2)}{8} | \widehat{u}|^2 + \frac{1}{4} \log (1+ |\xi |^2 )    |\widehat{u}|^2 \\
&\leq& 3 E_0(t, \xi ),
\end{eqnarray*}
which also  holds for $t>0$ and $\xi \in {\bf R}^n$.
\hfill
$\Box$

\begin{lem}\label{E-E}
$$ \frac{d}{dt}E(t, \xi) + \frac{\rho(\xi)}{2} E(t, \xi) \leq 0 , \quad t>0,\; \xi \in {\bf R}^n.$$
 \end{lem}
 
{\it Proof.}  \eqref{eq5}, \eqref{EFR} and Lemma \ref{E-Eo} imply that
\begin{eqnarray*}
\frac{\mathrm{d} }{\mathrm{d} t} E(t, \xi) + \frac{\rho(\xi)}{2} E(t, \xi)   &=& R(t, \xi) - F(t,\xi) + \frac{\rho(\xi)}{2} E(t, \xi) \\
&\leq &  R(t, \xi) - F(t,\xi) + \frac{3 \rho(\xi)}{2} E_0(t, \xi) \\
&=& \rho  (\xi ) |\widehat{u}_t|^2 - \log (1+ | \xi |^2) |\widehat{u}_t|  -  \rho(\xi) \log (1+ |\xi |^2 )  | \widehat{u}|^2 + \frac{3 \rho(\xi)}{4} |\widehat{u}_t|^2  \\ &+& \frac{3 \rho(\xi)}{4} \log (1+ | \xi |^2) |\widehat{u}|^2 \\
&=&  \left (  \frac{7 \rho(\xi)}{4} - \log (1+ | \xi |^2) \right ) |\widehat{u}_t|^2 - \frac{1}{4} \rho(\xi) \log (1+ |\xi |^2 )  | \widehat{u}|^2   \\
& \leq 0, &
\end{eqnarray*}
where we have just used the fact that 
$$  \frac{7 \rho(\xi)}{4} - \log (1+ | \xi |^2) = \left\{\begin{matrix}
\frac{ -1 }{8} \log (1+ | \xi |^2)  & \text { se } | \xi | \leq  \sqrt{e-1}, \\ 
\\
  \frac{7 }{8} - \log (1+ | \xi |^2) & \text { se } | \xi | > \sqrt{e-1},
\end{matrix}\right. $$ and the fact that $ \log (1+ | \xi |^2) \geq 1 $  for   $ | \xi | > \sqrt{e-1} $. Therefore,  $$\frac{7 }{8} - \log (1+ | \xi |^2) < \frac{ -1 }{8}   $$ 
 for  $ |\xi| > \sqrt{e-1}$.
 \hfill
$\Box$

Now we may note that  Lemma \ref{E-E} implies 
$$ E(t, \xi) \leq E(0, \xi) e^{ - \frac{\rho (\xi)}{2}t} .$$

Combining the last estimate with Lemma \ref{E-Eo} we arrive at the important proposition. 
$$ E_0(t, \xi) \leq 6 E_0(0, \xi) e^{ - \frac{\rho (\xi)}{2}t} ,$$
for all $t>0$ and $\xi \in {\bf R}^n$.

That is, using the definition of $E(t,\xi)$ we have obtained the important point wise estimates in the Fourier space.
\begin{pro}\label{density}\, It holds that
\begin{equation}\label{eq6}
  |\widehat{u}_t(t,\xi)|^2 + \log (1+ | \xi |^2)|\widehat{u}(t,\xi)|^2 \leq 6 \left ( |\widehat{u}_1(\xi)|^2 + \log (1+ | \xi |^2)|\widehat{u}_0(\xi)|^2  \right ) e^{ - \frac{\rho (\xi)}{2}t}, 
\end{equation}
for all $t>0$ and $\xi \in {\bf R}^n$, and
\begin{equation}\label{eq7}
|\widehat{u}(t,\xi)|^2 \leq 6 \left ( \frac{1}{\log (1+ | \xi |^2)} |\widehat{u}_1(\xi)|^2 + |\widehat{u}_0(\xi)|^2  \right ) e^{ - \frac{\rho (\xi)}{2}t},
\end{equation}
for all $t>0$ and $\xi \in {\bf R}^n, \xi \neq 0$.
\end{pro}


\subsection{Proof of Propositions \ref{L-energy} and \ref{L2-u}}

In this subsection, let us prove Propositions \ref{L-energy} and \ref{L2-u} by basing on the results of Proposition \ref{density}.\\ 
We first prove Proposition \ref{L-energy}. To begin with, applying the Plancherel theorem, and integrating the inequality \eqref{eq6} over ${\bf R}^n$ one has 
\begin{eqnarray}\label{AAA}
\left \| u_t (t, \cdot) \right \|_{L^2}^{2} &+& \left \| L^{1/2}u (t, \cdot) \right \|_{L^2}^{2} = \left \| \widehat{u}_t (t, \cdot) \right \|_{L^2}^{2} + \left \| \log^{1/2} (1+ | \cdot |^2)| \hat{u}(t, \cdot) \right \|_{L^2}^{2}\nonumber \\
 &=& \int _{R^n} \left(|\widehat{u}_t|^2 + \log (1+ | \xi |^2)|\widehat{u}|^2  \right)d\xi \nonumber\\ 
 &\leq& 6 \int _{R^n} \left ( |\widehat{u}_1|^2 + \log (1+ | \xi |^2)|\widehat{u}_0|^2  \right ) e^{ - \frac{\rho (\xi)}{2}t} d \xi \\
 &=& 6 \int _{| \xi | \leq 1} |\widehat{u}_1|^2  e^{ - \frac{\rho (\xi)}{2}t}  d \xi + 6 \int _{| \xi | \leq 1} \log (1+ | \xi |^2)|\widehat{u}_0|^2  e^{ - \frac{\rho (\xi)}{2}t} d \xi \nonumber \\
 &+& 6 \int _{ 1 < | \xi | \leq \sqrt{e-1}} |\widehat{u}_1|^2  e^{ - \frac{\rho (\xi)}{2}t}  d \xi + 6 \int _{ 1 < | \xi | \leq \sqrt{e-1}} \log (1+ | \xi |^2)|\widehat{u}_0|^2  e^{ - \frac{\rho (\xi)}{2}t} d \xi \nonumber\nonumber \\
 &+& 6 \int _{| \xi | > \sqrt{e-1} } |\widehat{u}_1|^2  e^{ - \frac{\rho (\xi)}{2}t}  d \xi + 6 \int _{| \xi | > \sqrt{e-1}} \log (1+ | \xi |^2)|\widehat{u}_0|^2  e^{ - \frac{\rho (\xi)}{2}t} d \xi \nonumber\\
 &=& 6(A_1 +A_2 +A_3)\nonumber,
\end{eqnarray}
with $A_i$ ($i=1,2,3$) according to the integrals on low, middle and high frequencies, respectively.

\subsection*{1) Estimate on the zone  $ | \xi | \leq 1 $}
On this zone we have  $\rho (\xi ) = \frac{1}{2} \log (1+|\xi |^2)$.

At this stage we assume that the initial data  $ u_0 , u_1 \in L^{1}({\bf R}^n)$ . Then   $ \widehat{u}_0 ,  \widehat{u}_1 \in L^{\infty }({\bf R}^n) $ and  
$$   \left \| \widehat{u}_0 \right \|_{\infty } \leq  \left \| u_0 \right \|_{1} \text{ } \text{ and }  \text{ }  \left \| \widehat{u}_1 \right \|_{\infty } \leq  \left \| u_1 \right \|_{1} . $$

Then, using the definition of $\rho(\xi)$  we may estimate the integrals on the low frequency region as follows.
\begin{eqnarray*}
A_1 &=&
 \int _{| \xi | \leq 1} |\widehat{u}_1|^2  e^{ - \frac{\log (1+ | \xi |^2)}{4}t}  d \xi + \int _{| \xi | \leq 1} \log (1+ | \xi |^2)|\widehat{u}_0|^2  e^{ - \frac{\log (1+ | \xi |^2)}{4}t} d \xi \\
&=&  \int _{| \xi | \leq 1} |\widehat{u}_1|^2  (1+ | \xi |^2)^{- \frac{t}{4}}  d \xi + \int _{| \xi | \leq 1} \log (1+ | \xi |^2)|\widehat{u}_0|^2 (1+ | \xi |^2)^{- \frac{t}{4}}    d \xi \\
&\leq &  \left \| \widehat{u}_1 \right \|_{\infty}^{2} \int _{| \xi | \leq 1}  (1+ | \xi |^2)^{- \frac{t}{4}}  d \xi +   \left \| \widehat{u}_0 \right \|_{\infty}^{2}  \int _{| \xi | \leq 1} \log (1+ | \xi |^2) (1+ | \xi |^2)^{- \frac{t}{4}}    d \xi \\
&\leq &  \left \| u_1 \right \|_{1}^{2} \int _{| \xi | \leq 1}  (1+ | \xi |^2)^{- \frac{t}{4}}  d \xi +   \left \| u_0 \right \|_{1}^{2}  \int _{| \xi | \leq 1} \log (1+ | \xi |^2) (1+ | \xi |^2)^{- \frac{t}{4}}    d \xi \\
&=& \left \| u_1 \right \|_{1}^{2} \omega _n \int _0^1  (1+ r^2)^{- \frac{t}{4}} r^{n-1} d r +   \left \| u_0 \right \|_{1}^{2} \omega _n \int _0^1 \log (1+r^2) (1+ r^2)^{- \frac{t}{4}}   r^{n-1} d r \\
& \leq & \left \| u_1 \right \|_{1}^{2} \omega _n \int _0^1  (1+ r^2)^{- \frac{t}{4}} r^{n-1} d r +   \left \| u_0 \right \|_{1}^{2} \omega _n \int _0^1  (1+ r^2)^{- \frac{t}{4}}   r^{n+1} d r \\
& \leq & \left \| u_1 \right \|_{1}^{2} \omega _n   I_{n-1}( t/4)   +   \left \| u_0 \right \|_{1}^{2} \omega _n  I_{n+1} (t/4) \\
& \leq & C_n \left(\left \| u_1 \right \|_{1}^{2}    t^{- \frac{n}{2}} +   \left \| u_0 \right \|_{1}^{2}  t^{- \frac{n+2}{2}} \right),\quad t \gg 1,
\end{eqnarray*}
because of the fact that $\log(1+r^2)  \leq r^2$ for all  $r \geq 0$, and Lemma \ref{general-p} with \eqref{2.1}, where $C_n$ is a positive constant depending only on $n$.

\subsection*{2)\,Estimate on the middle frequency zone $ 1 \leq | \xi | \leq \sqrt{e-1}$}

In this middle region we also have $ \rho (\xi ) = \frac{1}{2} \log (1+|\xi |^2) $ and we may estimate $\log(1+|\xi|^2)$ by 
$$ \log 2 \leq \log (1+|\xi |^2) \leq 1. $$
\noindent
Thus, one has for all $t>0$
\begin{eqnarray*}
A_2 &=& \int _{ 1 \leq | \xi | \leq \sqrt{e-1}} |\widehat{u}_1|^2  e^{ - \frac{\log (1+ | \xi |^2)
}{4}t}  d \xi + \int _{ 1 \leq | \xi | \leq \sqrt{e-1}} \log (1+ | \xi |^2)|\widehat{u}_0|^2  e^{ - \frac{\log (1+ | \xi |^2)
}{4}t} d \xi  \\
&\leq& \int _{ 1 \leq | \xi | \leq \sqrt{e-1}} |\widehat{u}_1|^2  e^{ - \frac{\log 2
}{4}t}  d \xi + \int _{ 1 \leq  | \xi | \leq \sqrt{e-1}} |\widehat{u}_0|^2  e^{ - \frac{\log 2}{4}t} d \xi \\
\vspace{0.3cm}
&\leq & 2^{-\frac{t}{4}} \left \| \widehat{u}_1 \right \|_{2}^{2} + 2^{-\frac{t}{4}} \left \| \widehat{u}_0 \right \|_{2}^{2} \; = \;2^{-\frac{t}{4}} \left( \left \| u_1 \right \|_{2}^{2} +  \left \|u_0 \right \|_{2}^{2}\right)\quad (t \geq 0).
\end{eqnarray*}

\subsection*{3) Estimate on the high frequency zone  $  | \xi | \geq \sqrt{e-1}$}
On this region we have $ \rho (\xi ) = \frac{1}{2} $. Thus we obtain the estimate 
\begin{eqnarray*}
A_3 &=& \int _{| \xi | \geq \sqrt{e-1} } |\widehat{u}_1|^2  e^{ - \frac{t}{4}}  d \xi + \int _{| \xi | \geq \sqrt{e-1}} \log (1+ | \xi |^2)|\widehat{u}_0|^2  e^{ - \frac{t}{4}} d \xi \\
&\leq & e^{ - \frac{t}{4}} \left \| \widehat{u}_1 \right \|_{2}^{2} + e^{ - \frac{t}{4}} \int _{{\bf R}^n} \log (1+ | \xi |^2)|\widehat{u}_0|^2   d \xi \\
& = & e^{ - \frac{t}{4}} \left(\left \| u_1 \right \|_{2}^{2} +  \left \| L^{\frac{1}{2}} u_0 \right \|_{2}^{2}\right) = 2e^{-\frac{t}{4}}E_u(0),  \quad t>0.
\end{eqnarray*}
Under these preparations obtained in 1), 2) and 3) above, one can prove Proposition \ref{L-energy}.\\

{\it Proof of Proposition \ref{L-energy}}\, By combining the estimates for $A_1, A_2, A_3$ with \eqref{AAA} the proof is now completed.
\hfill
$\Box$


The above proposition says that the total energy of the system decays as $t^{-n/2}$, that is

To estimate the $L^2$-norm of $u(t,x)$ we first observe that 
$$ \lim _{r \rightarrow 0} \frac{r^2}{\log (1+r^2)} = 1 .  $$
Thus,  there exists a small $\delta \in (0,1)$ such that 
$$ \frac{1}{2} \leq \frac{r^2}{\log (1+r^2)}\leq \frac{3}{2} $$
for $0< r \leq \delta$.

By integrating the inequality \eqref{eq7} on ${\bf R}^n$ and using the Plancherel theorem we obtain
 \begin{eqnarray}\label{BB}
 \left \| u(t,\cdot) \right \|_{2}^{2} & \leq & 6 \int _{{\bf R}^n} \left ( \frac{1}{\log (1+ | \xi |^2)} |\widehat{u}_1|^2 + |\widehat{u}_0|^2  \right ) e^{ - \frac{\rho (\xi)}{2}t} d \xi \nonumber\\
&=& 6 \int _{| \xi | \leq \delta } \left ( \frac{1}{\log (1+ | \xi |^2)} |\widehat{u}_1|^2 + |\widehat{u}_0|^2  \right ) e^{ - \frac{\rho (\xi)}{2}t} d \xi \\
 & +&  6 \int _{| \xi |  > \delta } \left ( \frac{1}{\log (1+ | \xi |^2)} |\widehat{u}_1|^2 + |\widehat{u}_0|^2  \right ) e^{ - \frac{\rho (\xi)}{2}t} d \xi \nonumber \\ 
 & =: & 6 (B_1+B_2).\nonumber
 \end{eqnarray}
 
Analogous to the estimates for the energy we may obtain exponential decay to the integral $B_2$ on the high frequency zone  $ | \xi| > \delta,$ 
that is, 
 $$B_2 \leq C \left( ||u_0||_2^2 + ||u_1||^2_2 \right) e^{-\frac{t}{4}}, \quad t>0.$$
On the low frequency region $ | \xi | \leq \delta $, by using Lemma \ref{general-p} together with \eqref{2.1} one has 
 \begin{eqnarray*}
 B_1 & =& \int _{| \xi | \leq \delta }  \frac{1}{\log (1+ | \xi |^2)} |\widehat{u}_1|^2 e^{ - \frac{\log (1+ | \xi |^2)}{4}t} d \xi + \int _{| \xi | \leq \delta }  |\widehat{u}_0|^2  e^{ - \frac{\log (1+ | \xi |^2)}{4}t} d \xi \\ 
& \leq &  \left \| u_1 \right \|_{1}^2 \int _{| \xi | \leq \delta }  \frac{1}{\log (1+ | \xi |^2)}  e^{ - \frac{\log (1+ | \xi |^2)}{4}t} d \xi + \left \| u_0 \right \|_{1}^2 \int _{| \xi | \leq \delta } (1+ | \xi |^2)^{- \frac{t}{4}}     d \xi \\ 
&\leq & \left \| u_1 \right \|_{1}^2 \omega _n \int _0^1  \frac{1}{\log (1+ r^2)}  (1+ r^2)^{- \frac{t}{4}}  r^{n-1} d r + \left \| u_0 \right \|_{1}^2 \omega _n \int _0^1(1+ r^2)^{- \frac{t}{4}} r^{n-1}     d r \\
&\leq & 
 \left \| u_1 \right \|_{1}^2 \frac{3\omega _n}{2} \int _0^1  (1+ r^2)^{- \frac{t}{4}}  r^{n-3} d r + \left \| u_0 \right \|_{1}^2 \omega _n \int _0^1(1+ r^2)^{- \frac{t}{4}} r^{n-1}     d r \\
 &=& \left \| u_1 \right \|_{1}^2 \frac{3\omega _n}{2} I_{ n-3}( t/4) + \left \| u_0 \right \|_{1}^2 \omega _n I_{ n-1}( t/4) 
\\
&\leq & C_n \left( \left \| u_1 \right \|_{1}^2  t^{- \frac{n-2}{2}} + \left \| u_0 \right \|_{1}^2  t^{- \frac{n}{2}}\right) , \quad  t \gg 1 
 \end{eqnarray*}
for $n>2$, where $C_n>0$ depends only on $n$.

{\it Proof of Proposition \ref{L2-u}}\, By combining estimates for $B_1, B_2$ with \eqref{BB}, we have just proved Proposition \ref{L2-u}.
\hfill
$\Box$


\section{Asymptotic profile of solutions}

To obtain an asymptotic profile we consider, without loss of generality, the case of initial amplitude $u_0 = 0$. Then, the corresponding Cauchy problem to problem \eqref{eqn}-\eqref{initial} in the Fourier space is given by 
\begin{align}\label{uhat}
&\hat{u}_{tt}(t,\xi) + \log(1+  |\xi|^2) \hat{u}_t (t,\xi) + \log(1+  |\xi|^2)\hat{u}(t,\xi) =0, \quad t > 0, \,\,\xi \in {\bf R}^{n},\,\\
&\hat{u}(0,\xi)=0, \quad
\hat{u_t}(0,\xi)=u_1(\xi), \quad \xi \in {\bf R}^{n}\nonumber.
\end{align}

The characteristics roots  $\lambda_+$ and  $\lambda_-$  of the characteristic polynomial 
$$\lambda^2+  \log(1+|\xi|^2)\lambda + \log(1+  |\xi|^2) =0, \quad  \xi \in {\bf R}^n$$
associated to the equation \eqref{uhat} are given by 

\begin{equation}\label{7.1}
\lambda_{\pm} = \dfrac{-\log(1+|\xi|^{2}) \pm i \sqrt{4\log(1+  |\xi|^2) - \log^2(1+|\xi|^{2}) }}{2},
\end{equation}
for $ \xi  \leq \sqrt{e^4-1} $. The solution formula can be expressed by  
\begin{equation}\label{expression}
\hat{u}(t,\xi)= \frac{\hat{u}_1}{b(\xi)} e^{-a(\xi)t}  \sin (b(\xi) t ) 
\end{equation}
for small frequency region such that $ | \xi | \leq \sqrt{e^4-1}$, where $ a(\xi)  $ and $ b(\xi) $ are the real and imaginary parts of the characteristics  roots, that is 
\begin{align}
a(\xi ) =  \frac{\log(1+|\xi|^{2})}{2} \text{ and } b (\xi) =  \frac{\sqrt{4\log(1+  |\xi|^2) - \log^2(1+|\xi|^{2}) }}{2}. 
\end{align}
\noindent
We note that $ a(\xi) $ and $ b(\xi)  $ are well defined for $ |\xi| \leq 1$. In fact, it is easy to see that 
$$ 4\log(1+  |\xi|^2) - \log^2(1+|\xi|^{2}) >0 $$ 
for $ 0 \leq | \xi | < \sqrt{e^4-1} $. 
\begin{rem}\label{remark4.1}{\rm It holds that
$$\sqrt{\log (1+|\xi|^2)} \leq 2 b(\xi) \leq 2 \sqrt{\log (1+|\xi|^2)}$$
for $ | \xi | \leq 1$.

To see this, we observe that 
\begin{eqnarray*}
b (\xi) &=&  \frac{\sqrt{4\log(1+  |\xi|^2) - \log^2(1+|\xi|^{2}) }}{2} \leq   \frac{\sqrt{4\log(1+  |\xi|^2) }}{2} \\ 
&=& \sqrt{\log(1+  |\xi|^2)},
\end{eqnarray*}
and for $ | \xi | \leq 1 \leq \sqrt{e^3 -1 }$, we have 

\begin{eqnarray*}
1 \leq | \xi |^2 +1 \leq e^3 &\Leftrightarrow & 0 \leq \log (1+\xi^2) \leq 3  \Leftrightarrow   \log ^2(1+\xi^2) - 3 \log (1+\xi^2) \leq 0 \\  &\Leftrightarrow &    \log (1+\xi^2) \leq  4 \log (1+\xi^2) - \log ^2(1+\xi^2),\\
\end{eqnarray*}
thus 
$$ \frac{\sqrt{\log (1+\xi^2) }}{2}  \leq b(\xi ) .   $$}
\end{rem}

In order to study an asymptotic profile of the solution to problem \eqref{eqn}--\eqref{initial} we consider a decomposition of the Fourier transformed initial data.

\begin{rem}\label{obs1}
{\rm Using  the  Fourier transform we can get a decomposition of the initial data $\hat{u}_1$  as follows: 
$$\hat{u}_1(\xi)=A(\xi)-iB(\xi)+P_1,\quad \xi \in {\bf R}^n, $$
where  $P_1, A, B$  are defined by  
\[P_{1} = \int_{{\bf R}^n}u_{1}(x) dx, \quad A(\xi)=\int_{{\bf R}^n}u_{1}(x)\big(1-\cos(\xi x) \big)dx, \quad B(\xi) =\int_{{\bf R}^n}u_{1}(x)\sin(\xi x)dx.\]}
\end{rem}

According to the above decomposition we can know the following lemma (see \cite{I-04}).
 
\begin{lem}\label{lema2.5.4}

Let $\kappa \in [0,1]$. For $u_{1} \in L^{1,\kappa}({\bf R}^{n})$ and $\xi \in {\bf R}^{n}$ it holds that 
$$|A(\xi)|\leq K|\xi|^\kappa\|u_{1}\|_{L^{1,\kappa}} \quad \text{ and } \quad |B(\xi)|\leq M|\xi|^\kappa\|u_{1}\|_{L^{1,\kappa}},$$
with positive constants $K$ and $M$ depending only on $n$.
\end{lem}

Let us capture a leading term of the solution based on \eqref{expression} and Remark \ref{obs1}. First, we apply the mean value theorem to get
\begin{equation}\label{Equacao4.5}
\sin (b(\xi)t) - \sin \left ( t \sqrt{\log(1+|\xi|^2)} \right ) = t \cos (\mu(\xi) t ) \left [ b(\xi) - \sqrt{\log(1+|\xi|^2)}\right],
\end{equation}
where $\mu (\xi) = \theta b(\xi) + (1- \theta)\sqrt{\log(1+|\xi|^2 }$ for some  $0<\theta <1$. By this reason, we can rewrite the solution formula \eqref{expression} as 
\begin{eqnarray*}
 \hat{u}(t,\xi) &=& \frac{A(\xi) - i B(\xi)}{b(\xi)} e^{-a(\xi)t}  \sin (b(\xi) t )  + \frac{P_1}{b(\xi)} e^{-a(\xi)t}  \sin \left ( t \sqrt{\log(1+|\xi|^2)} \right ) \\ &+& P_1 \frac{\left ( b(\xi) - \sqrt{\log(1+|\xi|^2)} \right )}{b(\xi)} e^{-a(\xi)t}  t \cos (\mu(\xi) t ). 
\end{eqnarray*}

Our goal in this section is to get decay estimates in time to the remainder therms defined above. To proceed with that we define the next $3$ functions
\begin{eqnarray*}
F_1(t, \xi) &=& \frac{ A(\xi) - i B(\xi)  }{b(\xi)} e^{-a(\xi)t}  \sin (b(\xi) t ),\\
F_2(t, \xi) &=& P_1 \frac{\left ( b(\xi) - \sqrt{\log(1+|\xi|^2)} \right )}{b(\xi)} e^{-a(\xi)t}  t \cos (\mu(\xi) t ),\\
F_3(t, \xi) &=&  \frac{P_1}{b(\xi)} e^{-a(\xi)t}  \sin \left ( t \sqrt{\log(1+|\xi|^2)} \right ) . 
\end{eqnarray*}

Then, we get 
$$ \hat u(t, \xi)  - F_3(t, \xi) = F_1(t, \xi) + F_2(t, \xi).  $$ 
We know that 
$$ \lim _{r\rightarrow 0} \frac{r^2}{\log (1+r^2)} =1,$$
so, there exists $0< \delta _1 <1 $ such that 
\[\frac{r^2}{\log (1+r^2)} <2\]
for all $0 <  r < \delta _1  $. By using this fact, Remark \ref{remark4.1}, Lemma \ref{lema2.5.4} with $ \kappa =1 $ and Lemma \ref{general-p} we obtain 
\begin{eqnarray*} 
\int _{| \xi | \leq \delta _1 } |F_1(t, \xi) |^2 d\xi &\leq& 4\int _{| \xi | \leq \delta _1 } \frac{ | A(\xi) - i B(\xi)|^2  }{\log (1+|\xi|^2)} e^{-2a(\xi)t}  \sin ^2(b(\xi) t ) d \xi \\
&\leq & 4\int _{| \xi | \leq \delta _1 } \frac{ ( |A(\xi)| + |B(\xi)|)^2  }{\log (1+|\xi|^2)} e^{-2a(\xi)t}   d \xi \\
&\leq & 4\int _{r \leq \delta _1 } \frac{ ( K + M)^2 |\xi|^{2 }  \| u_1 \|_{1,1}^2 }{\log (1+|\xi|^2)} (1+|\xi|^2)^{-t}   d \xi \\
&=& 4\omega _n ( K+M )^2 \| u_1 \|_{1,1}^2 \int _{r \leq \delta _1 } \frac{ r^{n+1} }{\log (1+r^2)} (1+r^2)^{-t}  d r \\ 
&\leq & 4\omega _n ( K+M )^2 \| u_1 \|_{1,1}^2 \int _{r \leq \delta _1} \frac{  r^2}{\log (1+r^2)} (1+r^2)^{-t} r^{n-1}  d r \\
&\leq &  8 \omega _n ( K+M )^2 \| u_1 \|_{1,1}^2 \int _{r \leq \delta _1 } (1+r^2)^{-t} r^{n-1}  d r \\
&\leq &  8 \omega _n ( K+M )^2 \| u_1 \|_{1,1}^2 \int _{r \leq 1 }  (1+r^2)^{-t} r^{n-1}  d r \\
&= &8 \omega _n \| u_1 \|_{1,1}^2 ( K+M )^2 I_{n-1} (t) \\
&\leq & C_{1,n} \| u_1 \|_{1,1}^2 t^{-\frac{n}{2}}.
\end{eqnarray*}
Now, we observe that for $ r := | \xi | < \sqrt{e^4-1} $, we have 
\begin{eqnarray*}
b(r) - \sqrt{\log (1+r^2)} &=& \frac{\sqrt{4 \log (1+r^2) - \log ^2 (1+r^2)}}{2} - \sqrt{\log(1+r^2)} \\
&=& \sqrt{\log(1+r^2)} \left ( \sqrt{1 - \frac{\log ^2(1+r^2)}{4 \log (1+r^2)}} - 1 \right ) \\
&=& \sqrt{\log(1+r^2)} \left ( \sqrt{1 - \frac{\log ^2(1+r^2)}{4 \log (1+r^2)}} - 1 \right ) \frac{\left ( \sqrt{1 - \frac{\log ^2(1+r^2)}{4 \log (1+r^2)}} + 1 \right )}{ \left ( \sqrt{1 - \frac{\log ^2(1+r^2)}{4 \log (1+r^2)}} + 1 \right )} \\
&=& -  \sqrt{\log(1+r^2)} \frac{\frac{\log ^2(1+r^2)}{4 \log (1+r^2)} }{ 1+  \sqrt{1 - \frac{\log ^2(1+r^2)}{4 \log (1+r^2)}}  }. 
\end{eqnarray*}
Thus
 \begin{eqnarray*}
\left | b(r) - \sqrt{\log (1+r^2)}  \right | &=& \left | \sqrt{\log(1+r^2)} \frac{\frac{\log ^2(1+r^2)}{4 \log (1+r^2)} }{ 1+  \sqrt{1 - \frac{\log ^2(1+r^2)}{4 \log (1+r^2)}}  }  \right |\\
&\leq & \sqrt{\log(1+r^2)} \frac{\log ^2(1+r^2)}{4 \log (1+r^2)}, \quad (\forall r < \sqrt{e^{4}-1}).
\end{eqnarray*}
By combining this fact  with Remark \ref{remark4.1}, we obtain 
$$ \frac{\left | b(r) - \sqrt{\log (1+r^2)}  \right | ^2}{|b(r)|^2} \leq \log ^2 (1+r^2)\quad (\forall r \leq 1).$$
Also, we know that
$$ \lim _{r \rightarrow 0^+} \frac{\log ^2(1+r^2)}{r^4 } = 1, $$
thus there exists $ 0< \delta < \delta_1 $ such that 
\[\frac{1}{2} \leq \frac{\log ^2(1+r^2)}{r^4 } \leq \frac{3}{2}\]
for all $0 \leq r \leq \delta$. Hence (see \eqref{2.1})

\begin{eqnarray*}
\int _{\xi \leq \delta } | F_2(t, \xi) | ^2 d \xi &=& \int _{\vert\xi\vert \leq \delta } | P_1|^2 \frac{\left | b(\xi) - \sqrt{\log(1+|\xi|^2)} \right |^2}{|b(\xi)|^2} e^{-2a(\xi)t}  t^2 \cos ^2 (\mu(\xi) t )  d \xi \\
&\leq & | P_1|^2 t^2  \int _{\vert\xi\vert \leq \delta } \frac{\left | b(\xi) - \sqrt{\log(1+|\xi|^2)} \right |^2}{|b(\xi)|^2} (1+ | \xi |^2)^{-t}    d \xi \\
&= & | P_1|^2 t^2 \omega _n \int _{0}^{\delta } \frac{\left | b(r) - \sqrt{\log(1+r^2)} \right |^2}{|b(r)|^2} (1+ r^2)^{-t}  r^{n-1}  d r \\
&\leq & | P_1|^2 t^2 \omega _n \int _{0}^{\delta }   (1+ r^2)^{-t} \frac{\log ^{2}(1+r^2)}{r^4}  r^{n+3}  d r \\
&\leq &  \frac{3\omega _n}{2}  | P_1|^2 t^2 \int _{0}^{\delta }   (1+ r^2)^{-t}  r^{n+3}  d r \\
&\leq &  \frac{3\omega _n}{2}  | P_1|^2 t^2 \int _{0}^{1 }   (1+ r^2)^{-t}  r^{n+3}  d r\\
&\leq & \frac{3\omega _n}{2}  | P_1|^2 t^2 I_{n+3}(t)\\
& \leq & C_{2,n} |P_1|^2 t ^{- \frac{n}{2}}. 
\end{eqnarray*}
\noindent
Here, note that 
$$ |\hat u(t, \xi)  - F_3(t, \xi) | ^2 = | F_1(t, \xi) + F_2(t, \xi)|^2 \leq 2 \left ( | F_1(t, \xi)|^2 + |F_2(t, \xi)|^2 \right ) .$$

Now, we define the following function
\begin{equation} \label{Profile}
\varphi(t, \xi) =  \frac{P_1}{\sqrt{\log(1+|\xi|^2)}} e^{-a(\xi)t}  \sin \left ( t \sqrt{\log(1+|\xi|^2)} \right ),
\end{equation}
which is equivalent to $F_3(t,\xi)$ according to the Remark \ref{remark4.1}.

Then we have the following result, which implies that the leading term of the Fourier transformed solution in the low frequency region is the very $\varphi(t,\xi)$. The result holds for all $n \geq 1$. 
\begin{theo}\label{teorema4.1}
Let $ n \geq 1 $  and let $ u_1 \in L^2({\bf R}^n) \cap L^{1,1}({\bf R}^n) $. Then, there exists $ 0 < \delta <1  $ such that 
$$ \int _{| \xi| \leq \delta } | \hat u(t, \xi)  - \varphi(t, \xi)| ^2 d \xi \leq  C_{1,n} \| u_1 \|_{1,1}^2 t^{-\frac{n}{2}} +  C_{2,n} |P_1 | ^2 t^{-\frac{n}{2}},$$
for $ t \gg 1$ with positive constants $C_{1,n} $ and $C_{2,n}$ depending only on $n \in \bf N$. 
\end{theo}

On the other hand, in the zone of high frequency $\{ |\xi| \geq \delta \}$ we have the following estimates.
\begin{theo}\label{teorema4.2}
Let $ n \geq 1 $  and let $ u_1 \in L^2({\bf R}^n) \cap L^{1}({\bf R}^n) $. Then,
$$ \int _{| \xi| \geq  \delta } | \hat u(t, \xi)  - \varphi(t, \xi)| ^2 d \xi \leq  C( \| u_1 \|^2 +  |P_1 | ^2 ) e^{-\eta t},$$
for $ t \gg 1$ with positive constant $C$ and $\eta$ depending on $n \in \bf N$. 
\end{theo}

{\it Proof.}\,It follows from Proposition \ref{density} that
\[\vert\hat{u}(t,\xi)\vert^{2} \leq 6\frac{\vert \hat{u}_{1}(\xi)\vert^{2}}{\log(1+\vert\xi\vert^{2})}e^{-\frac{t}{4}}\]
\begin{equation}\label{h-1}
\leq 6\vert \hat{u}_{1}(\xi)\vert^{2}e^{-\frac{t}{4}} \quad (\vert\xi\vert \geq \sqrt{e-1}).
\end{equation}
Additionally,   
\[\vert\varphi(t,\xi)\vert^{2} \leq \vert P_{1}\vert^{2}\frac{1}{\log(1+\vert\xi\vert^{2})}(1+\vert\xi\vert^{2})^{-t}\sin^{2}(t\sqrt{\log(1+\vert\xi\vert^{2})})\]
\begin{equation}\label{h-2}
\leq \vert P_{1}\vert^{2}(1+\vert\xi\vert^{2})^{-t}\quad  (\vert\xi\vert \geq \sqrt{e-1}).
\end{equation}
Thus, \eqref{2.2}, \eqref{h-1}, \eqref{h-2} and Lemma \ref{infit} imply
\[\int_{\vert\xi\vert \geq \sqrt{e-1}}\vert\hat{u}(t,\xi)- \varphi(t,\xi)\vert^{2}d\xi \leq 2\int_{\vert\xi\vert \geq \sqrt{e-1}}(\vert\hat{u}(t,\xi)\vert^{2} + \vert\varphi(t,\xi)\vert^{2})d\xi\]
\[\leq 12 e^{-\frac{t}{4}}\Vert u_{1}\Vert^{2} + 2\vert P_{1}\vert^{2}\omega_{n}\int_{1}^{\infty}(1+r^{2})^{-t}r^{n-1}dr\]
\[=  12 e^{-\frac{t}{4}}\Vert u_{1}\Vert^{2} +  2\vert P_{1}\vert^{2}\omega_{n}J_{n-1}\]
\begin{equation}\label{h-3}
\leq 12 e^{-\frac{t}{4}}\Vert u_{1}\Vert^{2} +  2\vert P_{1}\vert^{2}\omega_{n}\frac{2^{-t}}{t-1},\quad (t \gg 1).
\end{equation}
On the other hand, similarly to the derivation of \eqref{h-3}, if $\delta \leq \vert\xi\vert \leq \sqrt{e-1}$, then 
$$\log(1+\delta^{2}) \leq \log(1+r^{2}) \leq 1,$$
so that from Proposition \ref{density} one can get
\[\vert\hat{u}(t,\xi)\vert^{2} \leq 6\frac{\vert \hat{u}_{1}(\xi)\vert^{2}}{\log(1+\vert\xi\vert^{2})}e^{-\frac{\rho(\xi)}{2}t}\] 
\begin{equation}\label{h-4}
\leq  6\frac{\vert \hat{u}_{1}(\xi)\vert^{2}}{\log(1+\delta^{2})}e^{-\frac{\log(1+\delta^{2})}{4}t},
\end{equation}
and 
\[\vert\varphi(t,\xi)\vert^{2} \leq \frac{\vert P_{1}\vert^{2}}{\log(1+\vert\xi\vert^{2})}(1+\vert\xi\vert^{2})^{-t}\]
\begin{equation}\label{h-5}
\leq \frac{\vert P_{1}\vert^{2}}{\log(1+\delta^{2})}(1+\delta^{2})^{-t}\quad (\delta \leq \vert\xi\vert \leq \sqrt{e-1}).
\end{equation}
Finally, these estimates \eqref{h-3}, \eqref{h-4} and \eqref{h-5} in the high and middle frequency zones $\vert\xi\vert \geq \sqrt{e-1}$ and $\delta \leq \vert\xi\vert \leq \sqrt{e-1}$ imply the desired exponential decay estimates.
\hfill
$\Box$
\vspace{0.2cm}

The validity of Theorem \ref{main-theo} is a direct consequence of Theorems \ref{teorema4.1} ans \ref{teorema4.1}, and so we shall omit its detail.\\

\section{Optimal decay rate of $L^2$-norm for $n \geq 3$}

In this section, we investigate the precise rate of decay of the leading term \eqref{Profile} in $L^{2}$-sense as $t \to \infty$. The case of $n \geq 3$ is first treated.
\begin{pro}\label{proposition5.0}
Let $n \geq 3$. Then there exists $t_0 > 0$ such that for $t \geq t_0$ it holds that
$$ C_{1,n}t^{- \frac{n-2}{2}} \leq  \int _{{\bf R}^n} \frac{(1+|\xi |^2)^{-t} \sin ^2 (t\sqrt{ \log (1+|\xi |^2)})}{\log (1+|\xi |^2)} d\xi \leq C_{2,n}t^{- \frac{n-2}{2}},$$
where $ C_{1,n}$ and $ C_{2,n}$ are positive constants depending only on $n$. 
\end{pro}
{\it Proof.}\
We first observe that 
$$ \lim _{r \rightarrow 0} \frac{r^2}{\log (1+r^2)} = 1 .  $$
Then,  we can obtain  $ 0<\delta< 1  $   such that 
$$ \frac{1}{2} \leq \frac{r^2}{\log (1+r^2)}\leq \frac{3}{2} $$
for $0<r \leq \delta$. Thus, for $t \gg 1$ based on Lemmas \eqref{general-p} and \eqref{infit} one has
\begin{eqnarray*}
\int _{0}^{\delta} \frac{(1+r^2)^{-t} \sin ^2 (t\sqrt{ \log (1+r^2)})}{\log (1+r^2)} r^{n-1} dr &\leq & \int _{0}^{\delta} \frac{r^2}{\log (1+r^2)} (1+r^2)^{-t} r^{n-3} dr \\
& \leq & \frac{3}{2} \int _{0}^{\delta} (1+r^2)^{-t} r^{n-3} dr \\
& \leq & \int _{0}^{1} (1+r^2)^{-t} r^{n-3} dr \\
& \leq & C_{1,n} t^{- \frac{n-2}{2}},
\end{eqnarray*}

\begin{eqnarray*}
\int _{\delta}^{1} \frac{(1+r^2)^{-t} \sin ^2 (t\sqrt{ \log (1+r^2)})}{\log (1+r^2)} r^{n-1} dr &\leq  &\frac{1}{\log (1+\delta ^2)} \int _{\delta}^{1} (1+r^2)^{-t} \sin ^2 (t\sqrt{ \log (1+r^2)}) r^{n-1} dr \\
&\leq & \frac{1}{\log (1+\delta ^2)} \int _{\delta}^{1} (1+r^2)^{-t}  r^{n-1} dr \\
&\leq & \frac{1}{\log (1+\delta ^2)} \int _{0}^{1} (1+r^2)^{-t}  r^{n-1} dr \\
&\leq &  C_{2, n \delta } t^{-\frac{n}{2}} ,
\end{eqnarray*}

\begin{eqnarray*}
\int _{1}^{\infty} \frac{(1+r^2)^{-t} \sin ^2 (t\sqrt{ \log (1+r^2)})}{\log (1+r^2)} r^{n-1} dr &\leq & \frac{1}{\log 2} \int _{1}^{\infty} (1+r^2)^{-t} r^{n-1} \sin ^2 (t\sqrt{ \log (1+r^2)})  dr \\
 &\leq & \frac{1}{\log 2} \int _{1}^{\infty} (1+r^2)^{-t} r^{n-1}   dr \\
& \leq & C_{3,n}  \frac{2^{-t}}{t-1}.
\end{eqnarray*}
These three estimates above imply that there exists $t_0>0 $ such that  
\begin{equation}\label{largen}
\int _{0}^{\infty} \frac{(1+r^2)^{-t} \sin ^2 (t\sqrt{ \log (1+r^2)})}{\log (1+r^2)} r^{n-1} dr \leq  C_{4,n} t^{- \frac{n-2}{2}}
\end{equation}
for all $t\geq t_0$.

On the other hand, one notices the following computation such that
\begin{eqnarray*}
M(t) &:=& \int _{0}^{\infty} \frac{(1+r^2)^{-t} \sin ^2 (t\sqrt{ \log (1+r^2)})}{\log (1+r^2)} r^{n-1} dr \\ 
&=& \int _{0}^{\infty} \frac{(1+r^2)^{-t} (1+r^2) \sin ^2 (t\sqrt{ \log (1+r^2)})r^{n-2}\sqrt{t}r}{\sqrt{t}\sqrt{\log (1+r^2)} (1+r^2)\sqrt{\log (1+r^2)}}  dr \\
&\geq & \int _{0}^{\infty} \frac{(1+r^2)^{-t}  \sin ^2 (t\sqrt{ \log (1+r^2)})r^{n-2}\sqrt{t}r}{\sqrt{t}\sqrt{\log (1+r^2)} (1+r^2)\sqrt{\log (1+r^2)}}  dr . 
\end{eqnarray*}

And also, it is known that $ r^2 \geq \log (1+r^2) $, so that by using the change of variable  $y =\sqrt{t} \sqrt{\log (1+r^2) } $ we have 
\begin{eqnarray*}
M(t)&\geq & \int _{0}^{\infty} \frac{(1+r^2)^{-t}  \sin ^2 (t\sqrt{ \log (1+r^2)})(\log(1+r^2))^{\frac{n-2}{2}}\sqrt{t}r}{\sqrt{t}\sqrt{\log (1+r^2)} (1+r^2)\sqrt{\log (1+r^2)}}  dr  
\\
&=& \int _{0}^{\infty} \frac{e^{-y^2} \sin ^2 (\sqrt{t} y)y^{n-2}}{ t^{\frac{n-2}{2}}y }  dy \\
&=& \frac{t^{- \frac{n-2}{2}}}{2}\int _{0}^{\infty} e^{-y^2} y^{n-3} dy - \frac{t^{- \frac{n-2}{2}}}{2}\int _{0}^{\infty} e^{-y^2}y^{n-3} \cos (2\sqrt{t} y) dy \\
&=& \frac{t^{- \frac{n-2}{2}}}{2}(A_n - F_n(t)),
\end{eqnarray*}
where 
$$ A_n := \int _{0}^{\infty} e^{-y^2} y^{n-3} dy\,\, \text{and}\,\, F_n(t):=\int _{0}^{\infty} e^{-y^2}y^{n-3} \cos (2\sqrt{t} y) dy.  $$
Due to the fact $  e^{-y^2}y^{n-3} \in L^1({\bf R}) $ for $ n \geq 3$, we can apply the Riemann-Lebesgue Lemma to get 
$$ F_n (t) \rightarrow 0, \quad (t \to \infty) .$$
Then there exists $ t_1 >t_0  $ such that $ F_n(t) \leq \displaystyle{\frac{A_n}{2}}$ for all $ t \geq t_1  $, that is
$$ A_n - F_n(t) \geq  \frac{A_n}{2} \text{ for all } t \geq t_1. $$
Thus, one has
\begin{equation}\label{lagen2}
\int _{0}^{\infty} \frac{(1+r^2)^{-t} \sin ^2 (t\sqrt{ \log (1+r^2)})}{\log (1+r^2)} r^{n-1} dr \geq  \frac{A_n}{4} t^{- \frac{n-2}{2}}
\end{equation}
for $ t \geq t_1 $. 

Finally, the desired statement can be obtained by \eqref{largen} and \eqref{lagen2}.
\hfill
$\Box$


\section{Blow-up on infinite time for $n=1$ and  $n=2$}

In this section we study the optimal blowup rate in the sense of $L^{2}$-norm of the solution to problem (1.1)-(1.2).\\
We first derive the following proposition to the  case of dimension $n=1$.

\begin{pro}\label{proposition5.1}
There exists $ T >2  $  such that 
$$ \frac{(64+49\pi^2)t}{196 \pi ^2} \leq  \int _{\bf R} \frac{(1+|\xi |^2)^{-t} \sin ^2 (t\sqrt{ \log (1+|\xi |^2)})}{\log (1+|\xi |^2)} d \xi \leq 12 t . $$
for all $ t  \geq  T $. 
\end{pro}
{\it Proof.}\
We have 
$$ \frac{1}{2}\int _{\bf R} \frac{(1+|\xi |^2)^{-t} \sin ^2 (t\sqrt{ \log (1+|\xi |^2)})}{\log (1+|\xi |^2)} d \xi = \int _{0}^{\infty} \frac{(1+r^2)^{-t} \sin ^2 (t\sqrt{ \log (1+r^2)})}{\log (1+r^2)} d r .$$ 
Initially, we will obtain a lower bound for this integral. Set
$$ Q_l(t) : = \int _{0}^{\frac{1}{t}} \frac{(1+r^2)^{-t} \sin ^2 (t\sqrt{ \log (1+r^2)})}{\log (1+r^2)} dr,$$
$$ Q_h(t) : =  \int _{\frac{1}{t}}^{\infty} \frac{(1+r^2)^{-t} \sin ^2 (t\sqrt{ \log (1+r^2)})}{\log (1+r^2)} dr.$$
This implies
\[\frac{1}{2}{\cal I}_{1}(t) = Q_{l}(t) + Q_{h}(t).\]
By using the mean value theorem,  for $ 0 \leq r \leq \frac{1}{t} $  we obtain 
$$  \sin  (t\sqrt{ \log (1+r^2)}) \geq \frac{t}{2} \sqrt{ \log (1+r^2)} .$$

\begin{eqnarray*}
Q_l(t) &=& \int _{0}^{\frac{1}{t}} \frac{(1+r^2)^{-t} \sin ^2 (t\sqrt{ \log (1+r^2)})}{\log (1+r^2)} dr \\
&\geq & \frac{t^2}{4}   \int _{0}^{\frac{1}{t}} \frac{(1+r^2)^{-t} \log (1+r^2) }{\log (1+r^2)} dr \\ &=&
\frac{t^2}{4}  \int _{0}^{\frac{1}{t}}  (1+r^2)^{-t} dr \geq \frac{t^2}{4} \left ( 1+ \frac{1}{t^2} \right )^{-t}  \int _{0}^{\frac{1}{t}}  dr \\ 
&=& \frac{t}{4} \left ( 1+ \frac{1}{t^2} \right )^{-t}. 
\end{eqnarray*}

Now, since
$$ \lim_{t \rightarrow  \infty}  \left ( 1+ \frac{1}{t^2} \right )^{-t}  = 1 ,$$  
there exist a constant $t_1  \geq  1$ such that 
$$ \left ( 1+ \frac{1}{t^2} \right )^{-t}  \geq \frac{1}{2} \text{ for all } t > t_1 . $$
Thus for $ t  \geq  t_1 $, 
\begin{equation}\label{lowerql}
 Q_l (t) \geq  \frac{t}{8}.
\end{equation}

To deal with the integral $ Q_h(t) $, we consider 
$$ \nu _1 := \sqrt{e^{\frac{25 \pi ^2}{16 t^2}}-1} \text{ and } \nu _2 : = \sqrt{e^{\frac{49 \pi ^2}{16 t^2}}-1}. $$ 
Note that  for  $ \nu _1 \leq r \leq \nu _2  $  it holds that
$$ \left | \sin (t\sqrt{ \log (1+r^2)}) \right | \geq \frac{1}{\sqrt{2}} . $$
Then, we can estimate 
\begin{eqnarray*}
 Q_h(t)  &\geq &  \int _{\nu _1}^{\nu _2} \frac{(1+r^2)^{-t} \sin ^2 (t\sqrt{ \log (1+r^2)})}{\log (1+r^2)} dr\\
&\geq & \frac{1}{2} \int _{\nu _1}^{\nu _2} \frac{(1+r^2)^{-t} }{\log (1+r^2)} dr \\
&\geq &  \frac{8 t^2}{ 49 \pi ^2}  \int _{\nu _1}^{\nu _2} (1+r^2)^{-t} dr \\
&\geq & \frac{8 t^2}{ 49 \pi ^2} e^{\frac{-49 \pi^2}{16 t}} \int _{\nu _1}^{\nu _2} dr \\
&=&  \frac{8 t^2}{ 49 \pi ^2} e^{\frac{-49 \pi^2}{16 t}} \left ( \sqrt{e^{\frac{49 \pi ^2}{16 t^2}}-1} - \sqrt{e^{\frac{25 \pi ^2}{16 t^2}}-1} \right ). 
\end{eqnarray*}
Note that one knows the fact that
\[\lim_{t \to \infty}t\sqrt{e^{\frac{\gamma}{t^{2}}}-1} = \sqrt{\gamma}\quad (\gamma > 0).\]
Therefore, since one can get 
$$ \lim _{t \rightarrow \infty}  t e^{\frac{-49 \pi^2}{16 t}} \left ( \sqrt{e^{\frac{49 \pi ^2}{16 t^2}}-1} - \sqrt{e^{\frac{25 \pi ^2}{16 t^2}}-1} \right ) = \frac{\pi}{2}, $$
there exist $ t_2 \geq  1$ such that 
$$  t e^{\frac{-49 \pi^2}{16 t}} \left ( \sqrt{e^{\frac{49 \pi ^2}{16 t^2}}-1} - \sqrt{e^{\frac{25 \pi ^2}{16 t^2}}-1} \right ) \geq 1 \text{ for  } t  \geq  t_2. $$ 
Therefore, 
\begin{equation}\label{lowerqh}
Q_h(t)  \geq \frac{8 t}{ 49 \pi ^2} \text{ for } t \geq t_2.
\end{equation}
By adding \eqref{lowerql} and \eqref{lowerqh}, we conclude that
\begin{equation}\label{lowerestimateR}
\frac{1}{2}{\cal I}_{1}(t) =  \int _{\bf R} \frac{(1+|\xi |^2)^{-t} \sin ^2 (t\sqrt{ \log (1+|\xi |^2)})}{\log (1+|\xi |^2)} d \xi \geq  \frac{t}{4} + \frac{16 t}{ 49 \pi ^2},
\end{equation}
 for all $t  \geq  \text{max} \left \{ t_1, t_2 \right \}$.
 This estimate concludes the proof of lower bound of proposition.\\

\vspace{0.2cm}
 In order to obtain the upper bound, we separate the integral into three parts as follows:
$$ R_l(t) : = \int _{0}^{\frac{1}{t}} \frac{(1+r^2)^{-t} \sin ^2 (t\sqrt{ \log (1+r^2)})}{\log (1+r^2)} dr,$$
$$ R_m(t) : = \int _{\frac{1}{t}}^{\frac{1}{\sqrt{t}}} \frac{(1+r^2)^{-t} \sin ^2 (t\sqrt{ \log (1+r^2)})}{\log (1+r^2)} dr,$$
$$ R_h(t) : =  \int _{\frac{1}{\sqrt{t}}}^{\infty} \frac{(1+r^2)^{-t} \sin ^2 (t\sqrt{ \log (1+r^2)})}{\log (1+r^2)} dr.$$
This implies
\[\frac{1}{2}{\cal I}_{1}(t) = R_{l}(t) + R_{m}(t) + R_{h}(t).\]

Now, using the fact $\displaystyle{\frac{\vert\sin x\vert}{x}} \leq 1$ for all $x >0$, for $t>0$ one has
\begin{eqnarray}\label{R-low}
R_l(t)  
&\leq & \int _{0}^{\frac{1}{t}} (1+r^2)^{-t} t^2  dr \nonumber\\
&=& t^2 \int _{0}^{\frac{1}{t}}   (1+r^2)^{-t}   dr \\
&\leq & t^2 \int _{0}^{\frac{1}{t}}  dr \nonumber \\
&=& t . \nonumber
\end{eqnarray}

In order to estimate the middle part, we first observe that 
$$ \lim _{\sigma \rightarrow 0} \frac{\sigma }{\log (1 + \sigma )} = 1 . $$
So, there exists $ \delta _0 > 0 $ such that
$$  \frac{\sigma }{\log (1 + \sigma )} < 2 $$
for all $ 0 < \sigma < \delta _0 $. Therefore, if $\frac{1}{t} < r < \frac{1}{\sqrt{t}}$, then $ \frac{1}{t^2} < r^2 < \frac{1}{t}$ and for  $ t > \frac{1}{\delta_{0}}$, we have 
$$ \frac{ 1 }{\log (1 + r^2 )} < \frac{2}{r^2} .$$
Therewith,  using  integration by parts we can get 
\begin{eqnarray*}
R_m(t) &=& \int _{\frac{1}{t}}^{\frac{1}{\sqrt{t}}} \frac{(1+r^2)^{-t} \sin ^2 (t\sqrt{ \log (1+r^2)})}{\log (1+r^2)} dr \\
&\leq & 2 \int _{\frac{1}{t}}^{\frac{1}{\sqrt{t}}} \frac{(1+r^2)^{-t} }{r^2 } dr \\
&=& 2 t \left ( 1 + \frac{1}{t^2} \right )^{-t} - 2 \sqrt{t} \left ( 1 + \frac{1}{t} \right )^{-t}    - 4 t \int _{\frac{1}{t}}^{\frac{1}{\sqrt{t}}} (1+r^2 )^{-t-1} dr \\
&\leq & 2 t \left ( 1 + \frac{1}{t^2} \right )^{-t} . 
\end{eqnarray*}
\noindent
Since
$$\lim _{t \rightarrow \infty} \left ( 1 + \frac{1}{t^2} \right )^{-t} = 1, $$
there exists $ t_3  \geq 1 $ such that for all $ t  \geq t_3 $
$$  \left ( 1 + \frac{1}{t^2} \right )^{-t} \leq  2,$$
which implies
\begin{equation}\label{upperRm}
R_m(t) \leq 4t , \quad t  \geq t_3.
\end{equation}

Now, for $ t \geq 2 $ we estimate $R_{l}(t)$ as follows:
\begin{eqnarray*}
R_h(t) & =&  \int _{\frac{1}{\sqrt{t}}}^{\infty} \frac{(1+r^2)^{-t} \sin ^2 (t\sqrt{ \log (1+r^2)})}{\log (1+r^2)} dr \\
&\leq & \int _{\frac{1}{\sqrt{t}}}^{\infty} \frac{(1+r^2)^{-t} }{\log (1+r^2)} dr \\
&\leq &  \int _{\frac{1}{\sqrt{t}}}^{\infty} \frac{(1+r^2)^{-t+1}  }{ (1+r^2) \log (1+r^2)} dr \\
& \leq & \left ( 1+ \frac{1}{t} \right )^{-t+1} \frac{1}{\log \left ( 1+\frac{1}{t} \right )} \int _{\frac{1}{\sqrt{t}}}^{\infty} \frac{1  }{ 1+r^2 } dr \\
&=& \left ( 1+ \frac{1}{t} \right )^{-t+1} \frac{1}{\log \left ( 1+\frac{1}{t} \right )}   \left ( \frac{\pi}{2} - \tan ^{-1} (t^{-\frac{1}{2}}) \right ). 
\end{eqnarray*}
Due to the fact
$$ \lim _{t \rightarrow \infty}  \frac{1}{t}\left [ \left ( 1+ \frac{1}{t} \right )^{-t+1} \frac{1}{\log \left ( 1+\frac{1}{t} \right )}   \left ( \frac{\pi}{2} - \tan ^{-1} (t^{-\frac{1}{2}}) \right )  \right ]= \frac{\pi}{2e}, $$ 
there exist $ t_4  \geq  \max\{2, t_3\}$ such that 
$$  \left ( 1+ \frac{1}{t} \right )^{-t+1} \frac{1}{\log \left ( 1+\frac{1}{t} \right )}   \left ( \frac{\pi}{2} - \tan ^{-1} (t^{-\frac{1}{2}}) \right )  \leq  t $$
for all $t  \geq  t_4$, where one has just used the fact that
\[\lim_{t \to \infty}t\log(1+\frac{1}{t}) = 1.\]
Thus one has 
\begin{equation}\label{upperRh}
R_h(t) \leq t 
\end{equation}  
for all $t  \geq  t_4$. 

Finally, by adding \eqref{R-low}, \eqref{upperRm} and \eqref{upperRh} one can obtain the desired upper bound:
$${\cal I}_{1}(t) = \int _{0}^{\infty} \frac{(1+r^2)^{-t} \sin ^2 (t\sqrt{ \log (1+r^2)})}{\log (1+r^2)} d r \leq 6t$$
for all $t  \geq  t_4$.
\hfill
$\Box$
 \vspace{0.2cm}
  

Next we study the optimal blow-up order as $t \to \infty$ of ${\cal I}_{2}(t)$ given by 
$$ {\cal I}_{2}(t) =\int _{\bf R^2} \frac{(1+|\xi |^2)^{-t} \sin ^2 (t\sqrt{ \log (1+|\xi |^2)})}{\log (1+|\xi |^2)} d \xi  . $$
In order to do this we need the following elementary lemma.

\begin{lem}\label{Int-cos}
The  inequalities 
$$-1 \leq \int_2^{2\sqrt{t}} \frac{\cos y}{y}dy  \leq  1$$
hold for all $t>1 $.
\end{lem}

{\it Proof.}\,\, Using integration by parts we obtain for $t>1$

\begin{eqnarray*}
 \big|\int_2^{2\sqrt{t}} \frac{\cos y}{y}dy\big| &=& \big|\frac{1}{y}\sin y\Big|^{2\sqrt{t}}_{2} + 
\int^{2\sqrt{t}}_{2}\frac{1}{y^2}\sin y\; dy\big|\\
& \leq &  \frac{|\sin(2\sqrt{t})|}{2\sqrt{t}}+ \frac{|\sin 2|}{2}+ \int^{2\sqrt{t}}_{2}\frac{1}{y^2}\vert\sin y\vert dy\\
&\leq & \frac{1}{2\sqrt{t}} + \frac{1}{2} + \int^{2\sqrt{t}}_{2}\frac{1}{y^2}dy\\
&=& \frac{1}{2\sqrt{t}} + \frac{1}{2} + \frac{1}{2} - \frac{1}{2\sqrt{t}} = 1,
\end{eqnarray*}
which implies the desired estimate.
\hfill
$\Box$

\begin{rem}
{\rm We note that a more precise estimate than that in Lemma \ref{Int-cos} is 
$$-1 < \int_2^{2\sqrt{t}} \frac{\cos y}{y}dy  <  0, \quad t>1.$$
However, it is a little more difficult to prove.  For our propose in this paper, it is sufficient to use the rough estimate of Lemma \ref{Int-cos}. }
\end{rem}

\begin{pro}\label{proposition5.2}\, There exists $T>1$  such that
$$\frac{\pi}{4e} \log t  \leq   \int _{\bf R^2} \frac{(1+|\xi |^2)^{-t} \sin ^2 (t\sqrt{ \log (1+|\xi |^2)})}{\log (1+|\xi |^2)} d \xi \leq 6 \pi \log t $$
for all $ t \geq T$. 
\end{pro}

{\it Proof.}\,\,By considering the polar co-ordinate transform, we set 
$$\frac{1}{2\pi}{\cal I}_{2}(t) = \int _{0}^{\infty} \frac{(1+r^2)^{-t} \sin ^2 (t\sqrt{ \log (1+r^2)})}{\log (1+r^2)} r dr .$$

In order to obtain a lower bound for ${\cal I}_{2}(t) $, by using the change of variable $w=\sqrt{t\log(1+r^2)}$ and integration by parts, we observe that
\begin{eqnarray*}\label{lower-R2}
\frac{1}{2\pi}{\cal I}_{2}(t) 
&=& \int _{0}^{\infty} \frac{(1+r^2)^{-t} (1+r^2) \sin ^2 (t\sqrt{ \log (1+r^2)})\sqrt{t}r}{\sqrt{t}\sqrt{\log (1+r^2)} (1+r^2)\sqrt{\log (1+r^2)}}  dr \\
&\geq & \int _{0}^{\infty} \frac{(1+r^2)^{-t}  \sin ^2 (t\sqrt{ \log (1+r^2)})\sqrt{t}r}{\sqrt{t}\sqrt{\log (1+r^2)} (1+r^2)\sqrt{\log (1+r^2)}}  dr \\
&=& \int _{0}^{\infty} \frac{e^{-w^2} \sin ^2(\sqrt{ t}w)}{w}dw. \\
\end{eqnarray*}

Then one have 
\begin{eqnarray*}
\frac{1}{2\pi}{\cal I}_{2}(t) 
& \geq & 
 e^{-1} \int _{\frac{1}{\sqrt{t}}}^{1} \frac{ \sin ^2(\sqrt{ t}w)}{w}dw \\
&=& \frac{e^{-1}}{2} \int _{\frac{1}{\sqrt{t}}}^{1} \frac{ dw}{w} - \frac{e^{-1}}{2}  \int _{\frac{1}{\sqrt{t}}}^{1} \frac{ \cos (2 \sqrt{ t}w)}{w}dw \\
&= & \frac{e^{-1}}{4} \log t - \frac{e^{-1}}{2}  \int _{\frac{1}{\sqrt{t}}}^{1} \frac{ \cos (2 \sqrt{ t}w)}{w}dw 
\\
&=&  \frac{e^{-1}}{4} \log t - \frac{e^{-1}}{2}  \int _{2}^{2\sqrt{t}} \frac{ \cos y}{y}dy \\
&\geq& 
\frac{e^{-1}}{4} \log t - \frac{e^{-1}}{2}\\
&\geq & \frac{e^{-1}}{8} \log t , \quad  t\geq   e^4.
\end{eqnarray*}
The  penultimate inequality above is due to Lemma \ref{Int-cos}.
\vspace{0.2cm}

Thus, for $ t \gg 1 $, one has the optimal lower bound
\begin{equation}\label{2lower}
{\cal I}_{2}(t)=
\int _{0}^{\infty} \frac{(1+r^2)^{-t} \sin ^2 (t\sqrt{ \log (1+r^2)})}{\log (1+r^2)} r dr  \geq  \frac{\pi}{4 e} \log t.
\end{equation}
The estimate \eqref{2lower} implies the desired estimate from below of  Proposition  \ref{proposition5.2}.

\vspace{0.2cm}
Next, in order to get the upper bound for ${\cal I}_{2}(t)$ we set 
 $$ Q_l (t):=  \int _{0}^{\frac{1}{ t } } \frac{(1+r^2)^{-t} \sin ^2 (t\sqrt{ \log (1+r^2)})}{\log (1+r^2)} r dr , $$
$$ Q_m(t) :=  \int _{\frac{1}{ t } }^{\frac{1}{\sqrt{t}}} \frac{(1+r^2)^{-t} \sin ^2 (t\sqrt{ \log (1+r^2)})}{\log (1+r^2)} r dr,$$
$$  Q_h (t) :=  \int _{\frac{1}{ \sqrt{t} } }^{\infty} \frac{(1+r^2)^{-t} \sin ^2 (t\sqrt{ \log (1+r^2)})}{\log (1+r^2)} r dr . $$
This implies
\[\frac{1}{2\pi}{\cal I}_{2}(t) = Q_{l}(t) + Q_{m}(t) + Q_{h}(t).\]
\noindent
For $t>1$, we first have  
  \begin{eqnarray*}
  Q_l (t) &=&  \int _{0}^{\frac{1}{ t } } \frac{(1+r^2)^{-t} \sin ^2 (t\sqrt{ \log (1+r^2)})}{\log (1+r^2)} r dr \\
&\leq & \int _{0}^{\frac{1}{ t } } \frac{(1+r^2)^{-t} t^2 \log (1+r^2) }{\log (1+r^2)} r dr \\
& \leq & t^2  \int _{0}^{\frac{1}{ t } } (1+r^2)^{-t} r dr \\
& = & \frac{t^2 }{2 (t-1)}  \left [ 1- \left ( 1+\frac{1}{t^2} \right )^{1-t}\right].
 \end{eqnarray*}
Since 
 $$ \lim _{t\rightarrow  \infty} \frac{t^2 }{t-1}\left[ 1- \left (1+\frac{1}{t^2} \right)^{1-t} \right] = 1,$$ 
 there exists $t_2 \geq 1 $ such that 
 $$ \frac{t^2 }{2 (t-1)}  \left [ 1- \left ( 1+\frac{1}{t^2} \right )^{1-t} \right ] \leq 1$$
for all $t \geq  t_2$. 
 Therefore, it holds that
\begin{equation}\label{up-1}
Q_l (t) \leq 1
\end{equation}
 for $ t  \geq  t_2$. 
 \vspace{0.2cm}
 
Furthermore, for $t > 1$ one can get the estimate 
\begin{eqnarray*}
Q_m(t) &=&  \int _{\frac{1}{ t } }^{\frac{1}{\sqrt{t}}} \frac{(1+r^2)^{-t} \sin ^2 (t\sqrt{ \log (1+r^2)})}{\log (1+r^2)} r dr  \\
&\leq & \int _{\frac{1}{ t } }^{\frac{1}{\sqrt{t}}} \frac{r(1+r^2)^{-t} }{\log (1+r^2)}  dr  \\
&\leq & \int _{\frac{1}{ t } }^{\frac{1}{\sqrt{t}}} \frac{r(1+r^2)^{-1} }{\log (1+r^2)}  dr  \\
&=& \frac{1}{2} \left [  \log \left ( \log\left ( 1+\frac{1}{t} \right ) \right ) -  \log \left ( \log \left ( 1+\frac{1}{t^2} \right ) \right )\right]. 
\end{eqnarray*}   
Now, since we have 
$$ \lim_{t \rightarrow \infty } \frac{1}{\log t } \left [  \log \left ( \log\left ( 1+\frac{1}{t} \right ) \right ) -  \log \left ( \log \left ( 1+\frac{1}{t^2} \right ) \right ) \right ] = 1,$$
then there exists $ t_3 \geq t_2$ such that
$$ \frac{1}{2} \left [  \log \left ( \log\left ( 1+\frac{1}{t} \right ) \right ) -  \log \left ( \log \left ( 1+\frac{1}{t^2} \right ) \right ) \right ] \leq \log t$$
for all $t >t_3$, where one has just used the facts that
\[\lim_{\sigma \to +0}\frac{\log(\log(1+\sigma^{2}))}{\log\sigma} = 2,\]
\[\lim_{\sigma \to +0}\frac{\log(\log(1+\sigma))}{\log\sigma} = 1.\]
Therefore, one has just arrived at the estimate:
\begin{equation}\label{up-2}
Q_m(t) \leq \log t \quad (t  \geq t_3). 
\end{equation}

Similarly, for $t > 1$ it follows that
\begin{eqnarray*}
Q_h(t)&=&  \int _{\frac{1}{ \sqrt{t} } }^{\infty} \frac{(1+r^2)^{-t} \sin ^2 (t\sqrt{ \log (1+r^2)})}{\log (1+r^2)} r dr \\
&\leq &\int _{\frac{1}{ \sqrt{t} } }^{\infty} \frac{(1+r^2)^{-t} }{\log (1+r^2)} r dr \\
&\leq & \frac{1}{\log \left ( 1 + \frac{1}{t} \right )} \int _{\frac{1}{ \sqrt{t} } }^{\infty} (1+r^2)^{-t} r dr \\ 
&=& \frac{1}{2 (t-1) \log \left ( 1 + \frac{1}{t} \right )} \left (1 +\frac{1}{t} \right)^{1-t}.  
\end{eqnarray*}  
Since we see that
 $$ \lim _{t \rightarrow \infty}\frac{1}{(t-1) \log \left ( 1 + \frac{1}{t} \right )} \left ( 1 +\frac{1}{t} \right )^{1-t} = \frac{1}{e}, $$
there exists $ t_4 \geq t_3 > 1$ such that 
 $$ \frac{1}{2 (t-1) \log \left ( 1 + \frac{1}{t} \right )} \left ( 1 +\frac{1}{t} \right )^{1-t} \leq 1$$
for all $t  \geq t_4$. This implies 
 \begin{equation}\label{up-3}
 Q_h(t) \leq 1 \quad (t > t_4). 
 \end{equation}
By combining \eqref{up-1}, \eqref{up-2} and \eqref{up-3}, one can derive the crucial estimate
 \begin{equation}\label{2upper}
\frac{1}{2\pi}{\cal I}_{2}(t) = \int _{0}^{\infty} \frac{(1+r^2)^{-t} \sin ^2 (t\sqrt{ \log (1+r^2)})}{\log (1+r^2)} r dr \leq 3\log t
\end{equation}
for large $t  \geq  t_{4}$. 

The statement of Proposition \ref{proposition5.2} is now proved from \eqref{2lower} and \eqref{2upper}.

\hfill
$\Box$

\begin{rem}
{\rm The proof of Theorem \ref{main-theo2} is standard, and is a direct consequence of Theorem \ref{main-theo}, Propositions \ref{proposition5.0}, \ref{proposition5.1} and \ref{proposition5.2}. We omit its detail (see e.g., \cite{I-14}).}
\end{rem}
 

\par
\vspace{0.5cm}
\noindent{\em Acknowledgement.}
\smallskip
The work of the first author (R. C. CHAR\~AO) was partially supported by PRINT/CAPES - Process 88881.310536/2018-00 and the work of the third author (R. IKEHATA) was supported in part by Grant-in-Aid for Scientific Research (C)20K03682  of JSPS. 



\end{document}